\documentclass{amsart}
\pdfoutput=1
\usepackage{caption}
\usepackage{subcaption}
\usepackage{amssymb}
\usepackage{graphicx}
\usepackage{wrapfig}
\begin{document}

\author{G\'abor Fejes T\'oth}
\address{Alfr\'ed R\'enyi Institute of Mathematics,
Re\'altanoda u. 13-15., H-1053, Budapest, Hungary}
\email{gfejes@renyi.hu}

\author{L\'aszl\'o Fejes T\'oth}

\author{W{\l}odzimierz Kuperberg}
\address{Department of Mathematics \& Statistics, Auburn University, Auburn, AL36849-5310, USA}
\email{kuperwl@auburn.edu}

\title{Miscellaneous problems about packing and covering}
\thanks{The English translation of the book ``Lagerungen in der Ebene,
auf der Kugel und im Raum" by L\'aszl\'o Fejes T\'oth will be
published by Springer in the book series Grundlehren der
mathematischen Wissenschaften under the title
``Lagerungen---Arrangements in the Plane, on the Sphere and
in Space". Besides detailed notes to the original text the
English edition contains eight self-contained new chapters
surveying topics related to the subject of the book but not
contained in it. This is a preprint of one of the new chapters.}

\begin{abstract}
In this paper we discuss various special problems on packing and
covering. Among others we survey the problems and results concerning
finite arrangements, Minkowskian, saturated, compact, and totally
separable packings. We discuss shortest path problems and questions
about stability of packings.
\end{abstract}

\maketitle

\subsection{Arranging houses}

Suppose a large area is designated for a housing project in which the
minimum distance between the congruent rectangular outlines of the
houses is prescribed.  Which arrangement of the rectangles allows for
the greatest number of houses in the area? By (10,1), the problem is
reduced to the determination of the densest lattice packing of the parallel
domain of the rectangle. This problem was solved completely by {\sc{L.~Fejes
T\'oth}}~\cite{FTL67c} and {\sc{Florian}}~\cite{Florian68}. Let $a$
denote the length of the shorter side of the rectangle and $d$ the
prescribed minimum distance between the houses. There are three
types of optimal arrangement according as $a/d\le4-\sqrt{12}$,
$4-\sqrt{12}<a/d<2-\sqrt2$, or $2-\sqrt2\le{a/d}$. Let $H$ denote
the minimum area hexagon circumscribed about the parallel domain of
the rectangle. If $a/d\le4-\sqrt{12}$ then $H$ has bilateral symmetry
about a line parallel to the longer side of the rectangle (see Figure 1).
If $2-\sqrt2\le{a/d}$ then, besides a pair of sides parallel to the longer
sides of the rectangle, $H$ also has a pair of sides parallel to the shorter
sides of the rectangle (Figure 2). If $4-\sqrt{12}<a/d<2-\sqrt2$ then $H$
has neither bilateral symmetry nor sides parallel to the shorter sides of
the rectangle (Figure 3).

\medskip
\centerline {\immediate\pdfximage height3.38cm
{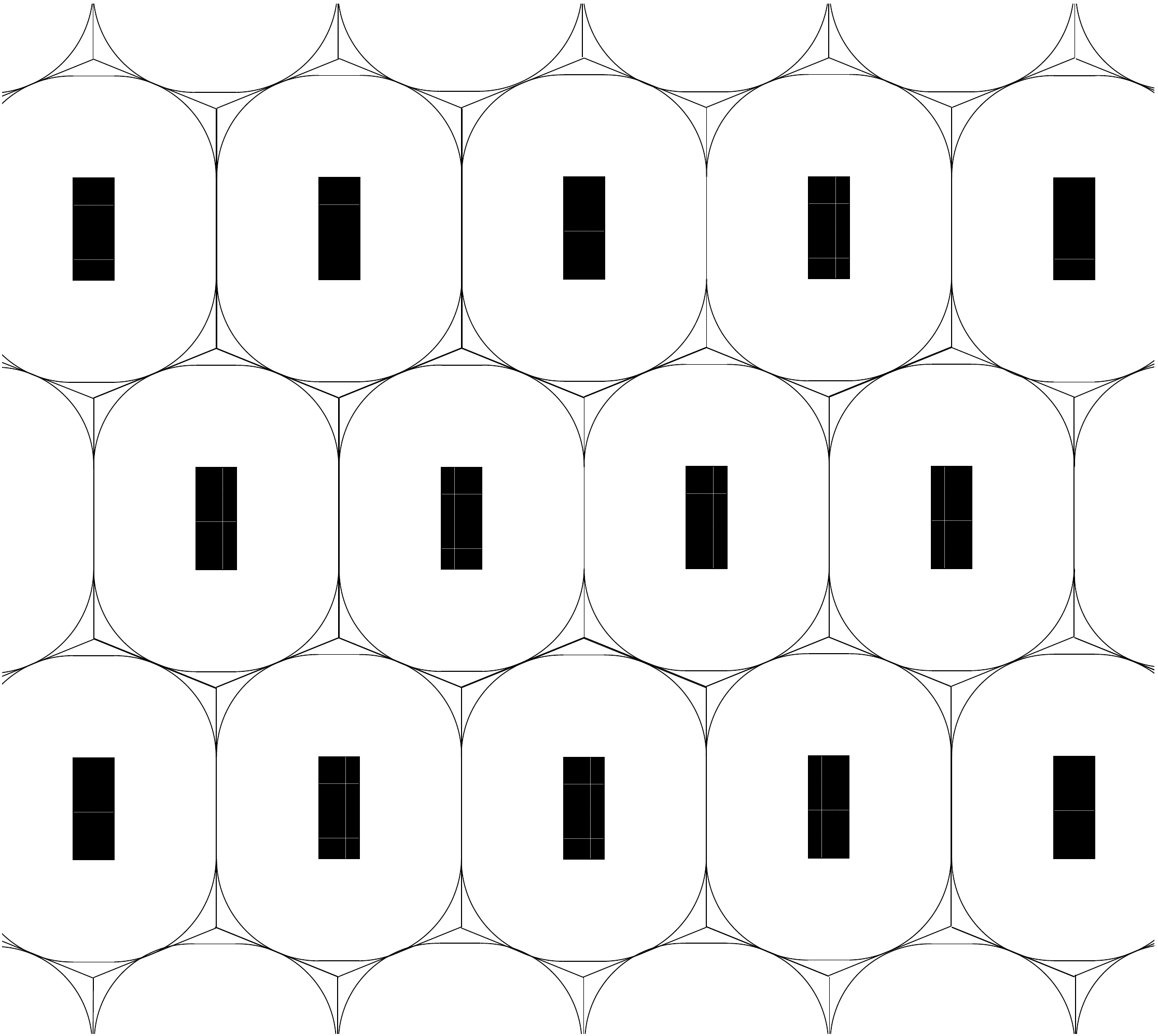}\pdfrefximage \pdflastximage\hskip0.6truecm
\immediate\pdfximage height3.38cm
{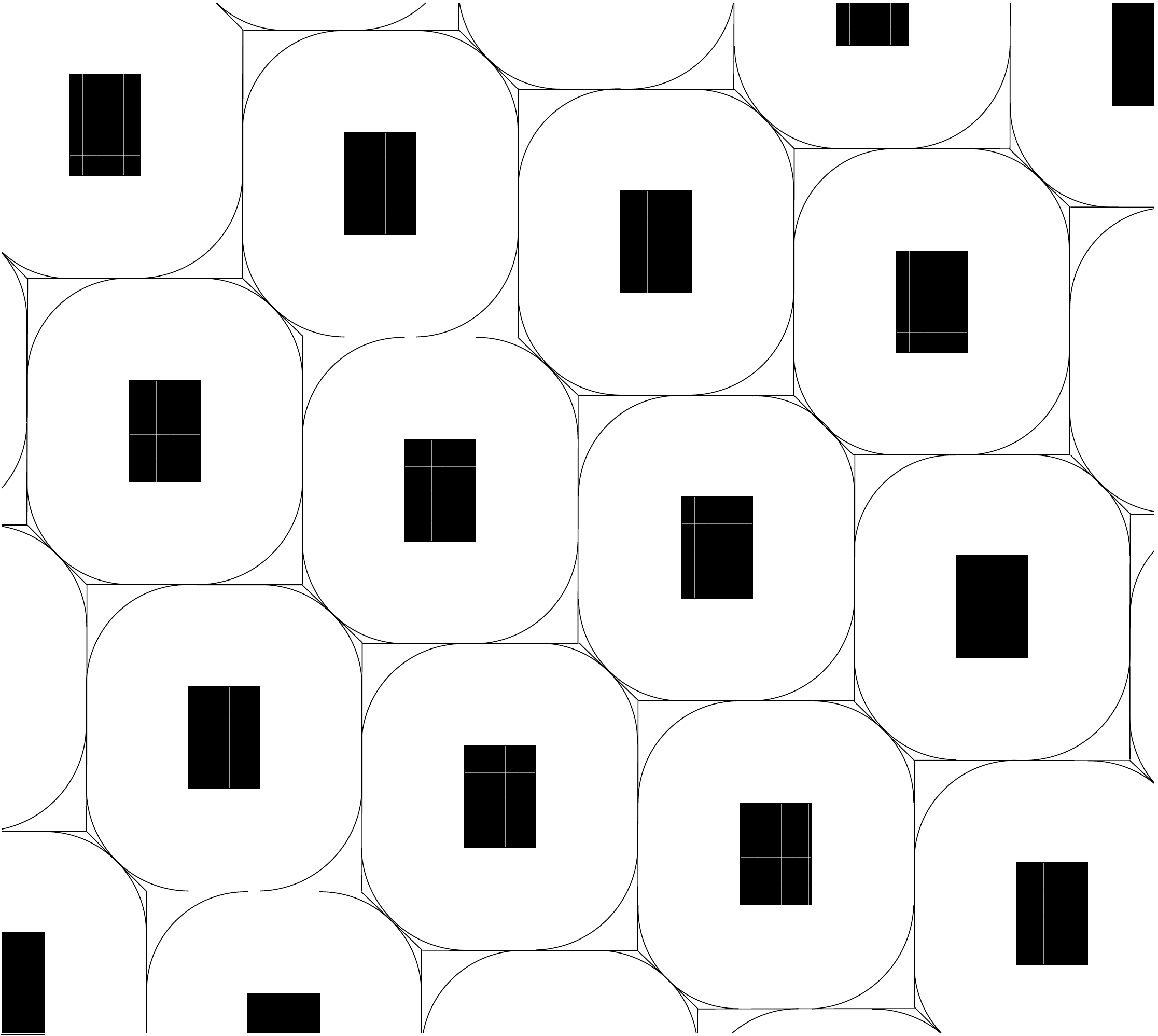}\pdfrefximage \pdflastximage\hskip0.6truecm
\immediate\pdfximage height3.38cm
{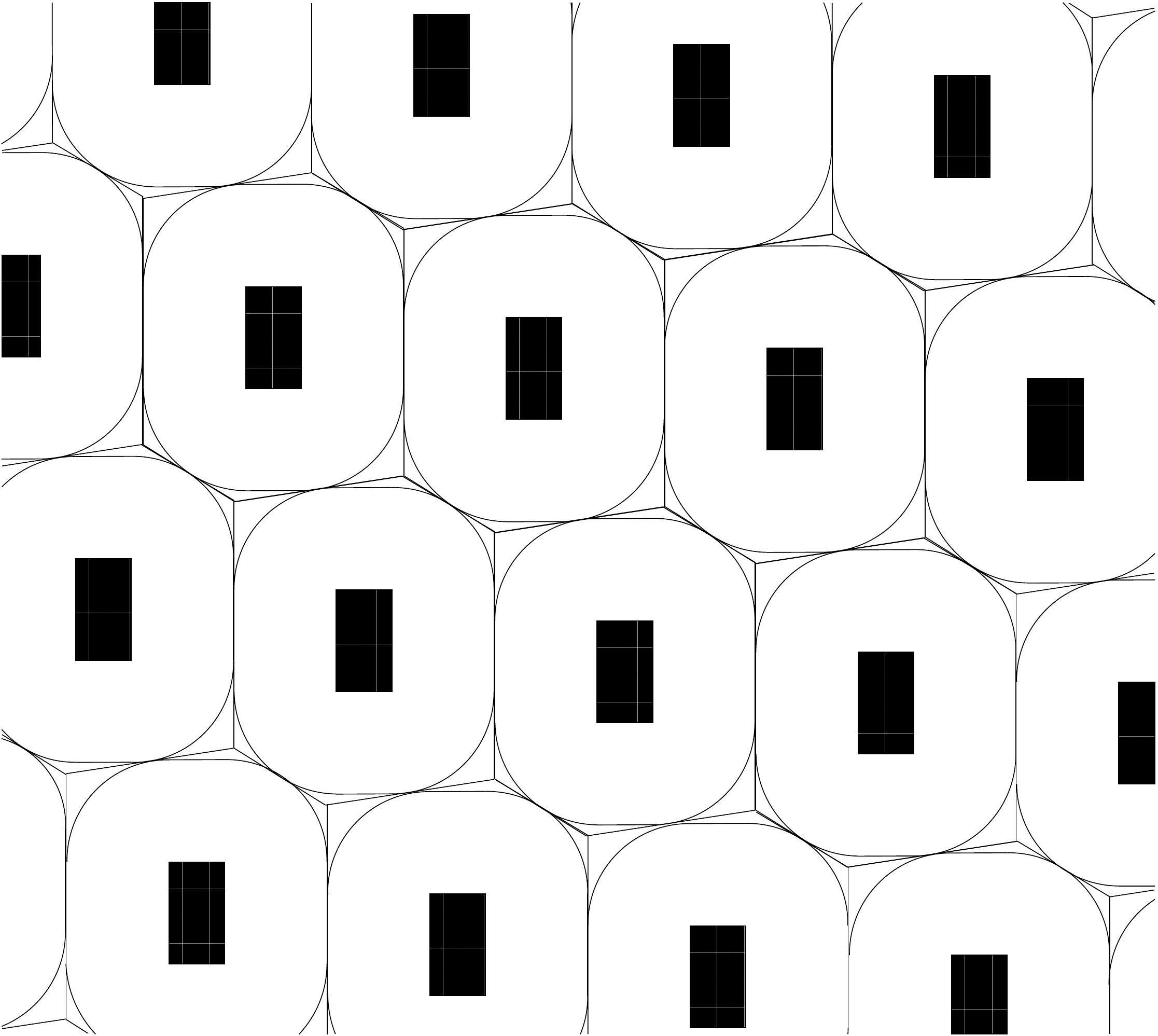}\pdfrefximage \pdflastximage}
\smallskip{\centerline{Figure~1 \hskip2.8truecm Figure~2 \hskip2.8truecm Figure~3}}
\medskip

Figure~4 shows a portion of an imaginary town. The black circles
represent cylindrical houses and the white ones landing platforms for helicopters.
To every black circle a white one is assigned, tangent to it. Naturally,
no black circle is allowed to overlap with another black circle or with
a white one. But the white circles can overlap with each other in part
or completely, so that several houses can share a common landing
platform. All white circles are mutually congruent and so are the black
ones, and their radii are prescribed. Under these conditions, the densest
packing of the black circles is to be determined. {\sc{Moln\'ar}}~\cite{Molnar64,Molnar66b}
and {\sc{Jucovi\v{c}}}~\cite{Jucovic70} obtained a series of very nice
arrangements of circles as the solution of this problem and of its variations.
The same problem on the sphere was studied by {\sc{Moln\'ar}}~\cite{Molnar75}

\medskip
\centerline {\immediate\pdfximage width4.55cm
{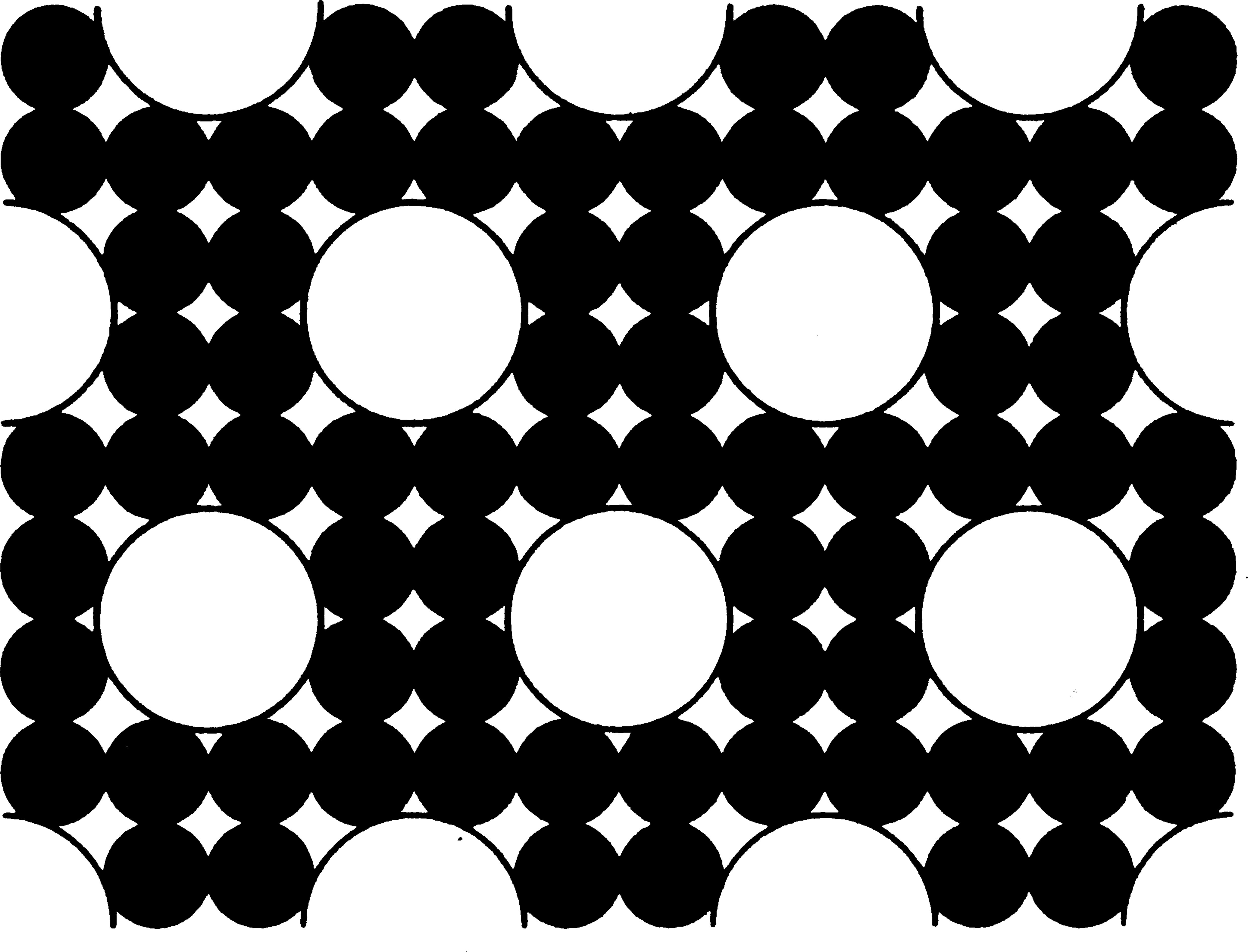}\pdfrefximage \pdflastximage}
\smallskip
{\centerline{Figure~4}}

\section{Packing barrels}

Let us imagine that in a large area we want to place as many
equal-sized barrels as possible, and so that each barrel can be
moved away from the area without disturbing the other barrels. Thus
the densest blocking-free packing of the plane with congruent circles
is desired, where a packing is {\it blocking-free} if within the part
of the plane not covered by the circles every circle can be moved
arbitrarily far from its original position. G.~Fejes T\'oth, who
raised this problem, stated the conjecture that in the best
packing the circles are arranged in double rows, separated by aisles
(Figure~5). The density of this packing equals
$\frac{\sqrt5-1}{2}\frac{\pi}{\sqrt{12}}$. {\sc{Heppes}}~\cite{Heppes67a} further
generalized the problem, asking only that the barrels be {\it
$r$-accessible}, meaning that a cellarer whose vertical shadow on the
floor is a circle of radius $r$ can freely come into contact with every
barrel. He conjectured that the densest $r$-accessible packing of
unit circles consists of appropriately placed double-rows. He gave
a density bound which brought him very
close to the solution of the problem. The correctness of the conjecture
was proved in a series of papers by {{\sc{G.~Blind}}
~\cite{Blind72,Blind76,Blind77}. A simpler proof was given by
\sc{G.~Blind}} and {\sc{R.~Blind}}~\cite{BlindBlind79}. {\sc{G.~Blind}}~\cite{Blind81} extended
the result to the sphere. The analogous problem in space was considered
by {\sc{G.~Blind}} and {\sc{R.~Blind}}~\cite{BlindBlind78}.

\medskip
\centerline {\immediate\pdfximage width5truecm
{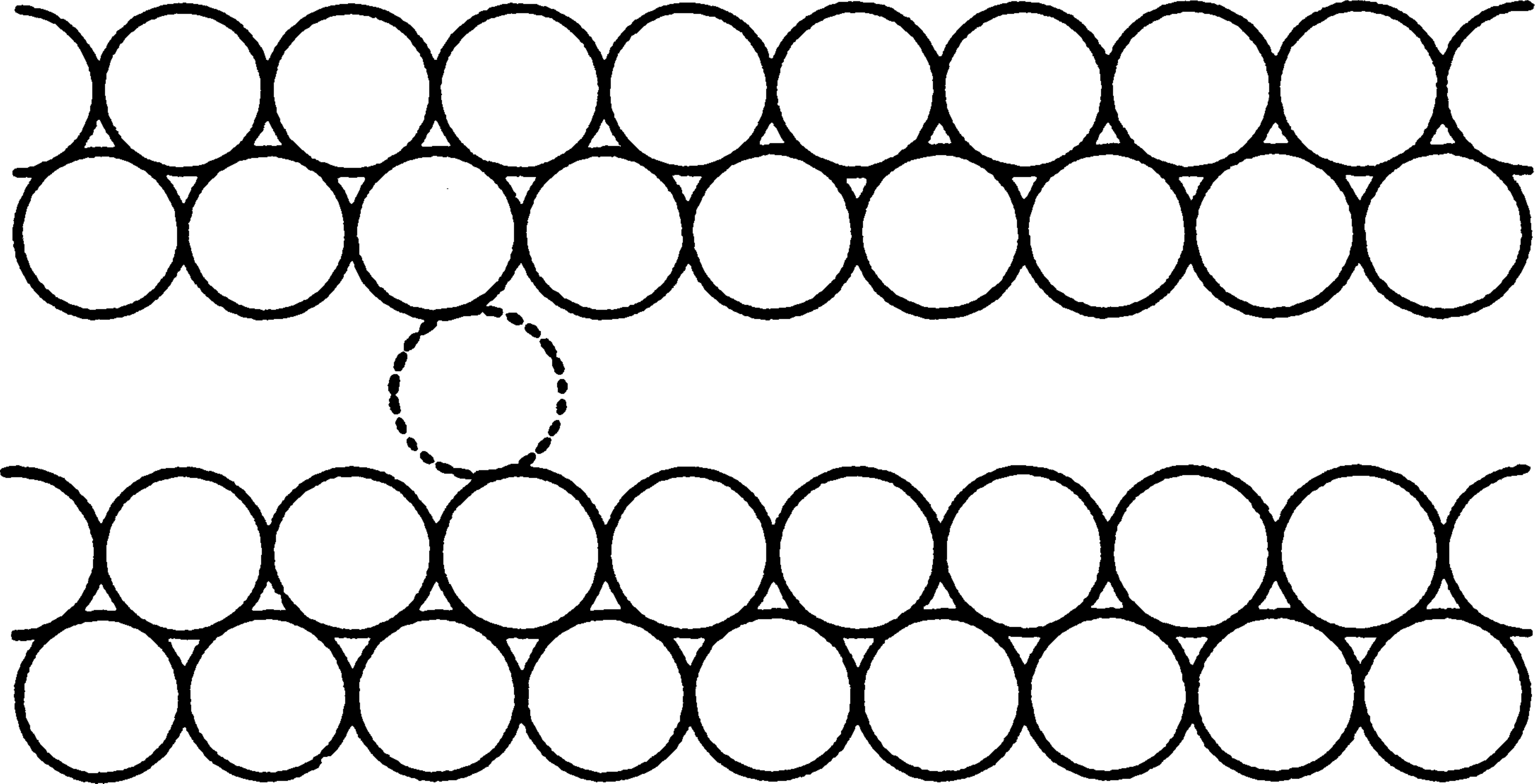}\pdfrefximage \pdflastximage}
\smallskip{\centerline{Figure~5}}

\section{Covering with a margin}

A covering of the plane with unit circles has {\it margin} $\mu\in[0,1]$ if the
removal of any one of the circles creates empty space small enough to be covered
by a circle of radius $1-\mu$. Obviously, a covering with a large margin cannot
be too thin. Hence, it is natural to ask: Determine the minimum
density of a covering with unit circles with margin $\mu$ and the covering that
attains this density. A covering with margin $0$ is a double covering, while
margin $1$ does not mean any additional restriction for the covering.
Thus, this problem connects the problems of thinnest covering and thinnest double
covering by unit circles.

{\sc{A.~Bezdek}} and {\sc{W.~Kuperberg}}~\cite{BezdekAKuperberg97}, who raised
the problem about thinnest covering with a margin, solved it restricted to the
special case when the original arrangement is lattice-like. For $0\le\mu\le\mu_1=0.56408\ldots$
the optimal arrangement is a triangular lattice, for $\mu_1\le\mu\le\mu_2=0.78608\ldots$
the solution is a square lattice, and for $\mu_2\le\mu\le1$ the optimal lattice remains
unchanged; it coincides with the thinnest double lattice of unit circles. Since the thinnest
double packing of unit circles is not lattice-like, the solution of the problem is
not lattice-like in general. However, it is conjectured that the solution is lattice-like
for sufficiently small values of $\mu$.

We note that for lattice arrangements the hole resulting by removing a circle
is symmetric about the center of the removed circle, thus the smallest circle
that can cover the hole will be centered here, as well. In this case
the problem can be interpreted as searching for the most
economical distribution of transmitting towers over a large area, all towers
having the same circular range, under the requirement that the region should be
covered even if due to a partial power loss the range of radius of one of the
towers is reduced by a factor of $1-\mu$.

{\sc Heppes}~\cite{Heppes02} considered a dual problem. A packing of unit circles has
{\it expendability} $\varepsilon>0$ if for every circle $C$ of the packing there is a
circle of radius $1+\varepsilon$ intersecting $C$ but not overlapping with any of the
other circles. Thus, removing $C$ and replacing it with the larger circle still creates
a packing. Heppes determined the densest lattice packing of unit circles with
expendability $\varepsilon$ for all $\varepsilon>0$.

\section{Finite packing and covering in 2 dimensions}

Nice, particular problems arise when trying to pack a given
number of congruent circles of maximum radius in a specific region
or to cover the region with congruent circles of minimum radius. Most often,
the chosen container is the square, the circle, or the equilateral
triangle.  The extensive literature on this subject consists mainly of articles
treating a single, specific case of the general problem, too numerous to be
listed here.  Among them are several algorithmic results that present
some good arrangements, however not confirmed to be optimal.  A major
goal in this field is to find algorithms that give or approximate the optimal
solutions.  Thus far, the only algorithm of this type was constructed by
{\sc{Peikert}}, {\sc{W\"urtz}}, {\sc{Monagan}} and {\sc{de Groot}} \cite{Peikert+}
(see also {\sc{Peikert}} \cite{Peikert}) for finding the densest packing of congruent
circles in a square. Figure~6 illustrates the
optimal packing of up to 20 congruent circles in a square container.

\medskip
\centerline {\immediate\pdfximage width9.5cm
{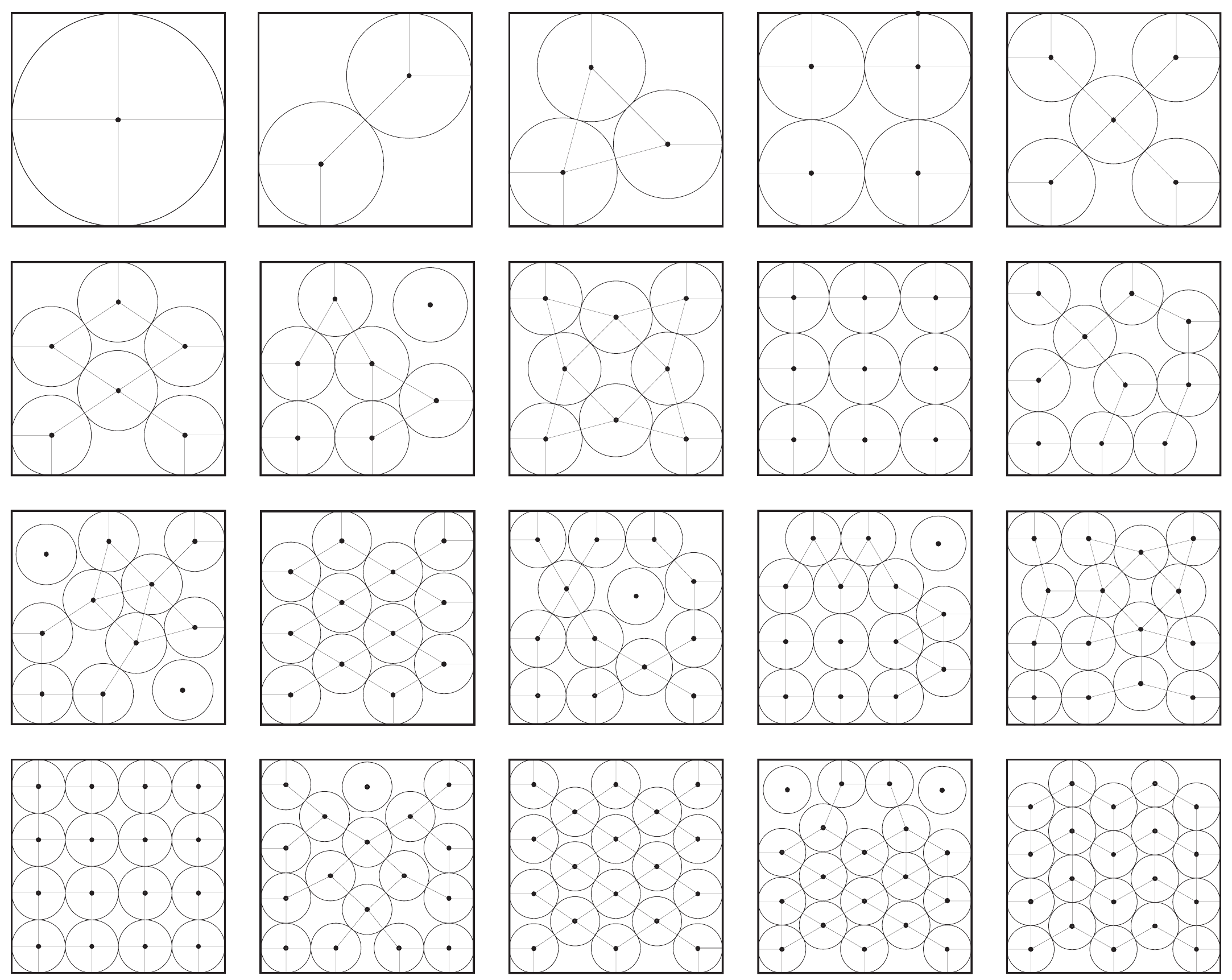}\pdfrefximage \pdflastximage}
\smallskip{\centerline{Figure~6}}
\medskip

Consider a densest packing with a large number of translates of a convex
disk in a large container, or under some other constraints that force
finiteness. While it is expected that such packing should be close to a cluster
taken from the densest packing of the whole plane or space, it seldom is
identical with such a cluster. {\sc{Sch\"{u}rmann}} \cite{Schurmann02a} proved that
the solutions to several finite packing problems are non-lattice if
the number of the translates is sufficiently large. In particular, he proved
this for the packings of circles that minimizes the diameter of their union,
thereby confirming a conjecture of {\sc Erd\H{o}s}.

For further literature on packings in bounded containers, see {\sc Szab\'{o}} {\it et
al.} \cite{Szabo+} and {\sc Melissen} \cite{Melissen97}.

\section{Finite arrangements in higher dimensions}

There are several results about packing congruent balls in various containers.
{\sc Schaer} \cite{Schaer66a,Schaer66b,Schaer66c,Schaer94} considered the
problem of densest packing of $k$ congruent balls in a cube, and solved it for
$k\le 10$. The notoriously difficult case of $14$ balls in a cube was settled
by {\sc Jo\'{o}s} in \cite{Joos09a} by proving that if $14$ points are placed
in the unit cube, then two of the points are no more than $\sqrt2/2$ away from
each other. The results of Schaer and Jo\'{o}s confirm some of the conjectures
stated by {\sc{Goldberg}} \cite{Goldberg71}.

{\sc{Golser}} \cite{Golser} studied the problem of packing $k$ congruent balls
in the regular octahedron and solved it for $k\le7$. {\sc{B\"{o}r\"{o}czky Jr.}}
and {\sc{Wintsche}} \cite{BoroczkyJrWintsche00} generalized Golser's result to
higher dimensions. They proved that the maximum radius of $k\le 2n+1$
congruent balls packed in the regular $n$-dimensional cross-polytope does not
depend on $n$, and for $4\le k\le 2n$ the radius is constant. {\sc{K.~Bezdek}}
\cite{BezdekK87a} solved the problem of packing $k$ congruent balls in a regular
tetrahedron, for $k=5, 8, 9$ and $10$. {\sc{W.~Kuperberg}} \cite{Kuperberg07}
considered the problem of maximum radius of $k\le 2n+2$ $n$-dimensional
congruent balls packed in a spherical container. While for $k\le n+1$ and
$k=2n+2$, the optimal configurations of balls are unique (see {\sc{Davenport}}
and {\sc{Haj\'{o}s}} \cite{DavenportHajos} and {\sc{Rankin}} \cite{Rankin55}),
W.~Kuperberg describes the structure of the non-unique configurations for
$n+2\le k\le 2n+1$, in which the radius remains constant. An alternate
characterization of these configurations was given by {\sc{Musin}}
\cite{Musin19}.

The problem of thinnest covering of the cube with $k$ congruent balls was
solved by {\sc{G. Kuperberg}} and {\sc{W. Kuperberg}} for $k=2,\,3,\,4$ and $8$
and by {\sc Jo\'os} \cite{Joos14a,Joos14b} for $k=5$ and 6. {\sc Jo\'{o}s}
\cite{Joos08b,Joos09b,Joos18} proved that the minimum radius of $8$ congruent
balls that can cover the unit cube is $\sqrt{5/12}$ in $4$ dimensions and
$\sqrt{2/3}$ in $5$ dimensions. The problem of covering the $n$-dimensional
cross-polytope with $k$ congruent balls of minimum radius was studied by
{\sc{B\"{o}r\"{o}czky, Jr.}}, {\sc{F\'{a}bi\'{a}n}} and {\sc{Wintsche}}
\cite{BoroczkyJrFabianWintsche} who found the solution for $k=2,\ n$, and $2n$.
Remarkably, the solution of the cases $k=2$ and $k=n$ is substantially different
for $n=3$ and $n\neq3$.

Another finite packing problem asks to arrange $k$ non-overlapping unit balls
so that the convex hull of their union is of minimum volume.  {\sc L.~Fejes
T\'{o}th} \cite{FTL75b} conjectured that in dimension $n\ge5$ the balls' centers
should be collinear, so that the convex hull of the union of the balls forms a
``sausage-like'' solid of length $2k$ (see Figure~7).

\medskip
\centerline {\immediate\pdfximage width8cm
{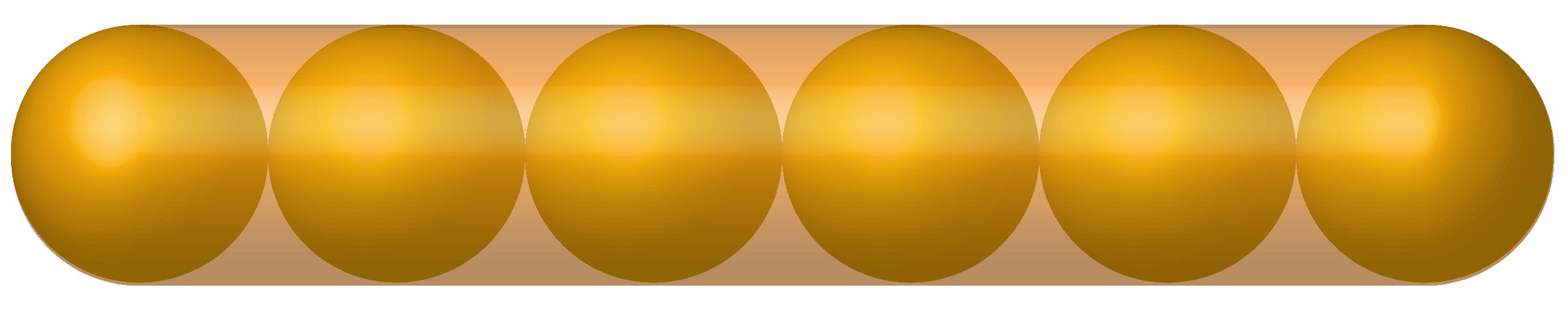}\pdfrefximage \pdflastximage}
\smallskip{\centerline{Figure~7}
\medskip

The conjecture, known as the {\it{sausage conjecture}}, attracted great interest and
generated intensive research on finite packing and covering (see
{\sc{Gritzmann}} and {\sc{Wills}} \cite{GritzmannWills} and {\sc{B\"or\"oczky Jr.}}
\cite{BoroczkyJr04}). The conjecture was verified for $n\ge13,387$ by
{\sc Betke}, {\sc Henk} and {\sc Wills} \cite{BetkeHenkWills}. {\sc Betke}
and {\sc Henk} \cite{BetkeHenk98} lowered the bound for the dimension to $n\ge42$.

\section{Slab, cylinder, torus}

What is the densest packing of unit circles in a strip of width $2\ge{w}$?
The answer is trivial for $2\le{w}\le2+\sqrt3$ but it becomes difficult
for ${w}\ge2+\sqrt3$. Extending a result by {\sc Kert\'esz} \cite{Kertesz82}, who
gave the solution for $2+\sqrt3\le w\le 2+2\sqrt2$, {\sc F\"uredi} \cite{Furedi91}
solved the problem for $2+\sqrt3\le w\le 2+2\sqrt3$. {\sc{Moln\'ar}}
conjectured that the maximum density is
$$\frac{(n+1)(n+2)\pi}{2w(n+1)+2\sqrt{1-(w-\sqrt3)^2}},\quad
n=\left\lfloor\frac{w-2}{\sqrt3}\right\rfloor,$$
and observed that for $w=2+n\sqrt3$, $n=1,2,\ldots$ his conjecture follows
from the theorem of {\sc{Groemer}} \cite{Groemer60a} mentioned on page XX.

{\sc Horv\'ath} and {\sc Moln\'ar} \cite{HorvathMolnar67} studied the problem
of densest packing of unit balls in a slab of space bounded by a pair of
parallel planes. They showed, among other things, that such packings
consisting of two hexagonal or square layers of balls are extremal
in a slab of the corresponding width. {\sc Moln\'ar} \cite{Molnar78} extended
this result by finding the densest packing in a slab of every width
between $2$ and $2+\sqrt2$. \sc{Horv\'ath}} \cite{Horvath74} investigated
the problem of densest packing of unit balls in a 4-dimensional slab of
width $2<w\le2+\sqrt2$.

Packing and covering with congruent circles on the surface of the infinite circular
cylinder was considered by {\sc L.~Fejes T\'{o}th} \cite{FTL62}, {\sc Bleicher} and
{\sc L.~Fejes T\'{o}th} \cite{BleicherFTL64} and {\sc{Mughal}} and {\sc{Weaire}}
\cite{MughalWeaire}. {\sc{L.~Fejes T\'oth}} \cite{FTL73e} gave an upper bound for
the number of points with given minimal distance on the surface of polyhedra.
The papers by {\sc{Dickinson, Guillot, Keaton}} and {\sc{Xhumari}}
\cite{DickinsonGuillotKeatonXhumari11a,DickinsonGuillotKeatonXhumari11b},
{\sc{Connelly}} and {\sc{Dickinson}} \cite{ConnellyDickinson},
{\sc{Connelly, Shen,}} and {\sc{Smith}} \cite{ConnellyShenSmith},
{\sc{Connelly, Funkhouser, V. Kuperberg}} and {\sc{Solomonides}}
\cite{ConnellyFunkhouserKuperbergSolomonides}, {\sc{Heppes}}
\cite{Heppes99}, {\sc{Musin}}, {\sc{Nikitenko}} \cite{MusinNikitenko},
{\sc{Przeworski}} \cite{Przeworski06b} and {\sc{Brandt, Dickinson, Ellsworth, Kenkel}}
and {\sc{Smith}} \cite{BrandtDickinsonEllsworthKenkelSmith} investigate
packings of circles on the torus. The dual problem of the thinnest covering of
the torus by congruent circles was treated by {\sc{Jo\'os}} \cite{Joos19}.
{\sc{Jo\'os}} and {\sc{Nagy}} \cite{JoosNagy} determined the smallest upper
bound for the radius of of $k\le4$ congruent balls in the 3-dimensional cubical flat torus.

\section{Close packings and loose coverings}

{\sc{L.~Fejes T\'oth}}~\cite{FTL76} defined another measure of efficiency,
alternate to density. He considered the supremum of the radii of the circles
that can be placed in the complement of the packing. The smaller that number
is, the more close, or efficient, is the packing. The {\it{closeness}} of the
packing is measured by the inverse of this supremum, and a {\it close} packing
is one with largest possible closeness. {\it Looseness of a covering} is similarly
determined by the inverse of the supremum of the radii of circles that can be
placed in the intersection of two members of the covering, and a {\it loose}
covering is one with largest possible looseness. {\sc{L.~Fejes T\'oth}}
~\cite{FTL78a} proved that the closeness
of a packing by translates of a convex disk $K$ cannot exceed the
closeness of the closest lattice packing of $K$. He also remarked that
this remains true if positively-homothetic copies of another convex disk
instead of a circle are used to measure closeness. This is a result analogous
to the corresponding theorem about density. On the other hand, he produced
a centrally symmetric convex disk and a packing consisting of translates of
the disk and a rotated copy of it, with closeness greater than that
of the closest lattice packing. An alternative example was constructed
by {\sc{A.~Bezdek}}~\cite{BezdekA80}.

{\sc{Linhart}}~\cite{Linhart78} observed that the natural way of measuring
closeness of a packing and looseness of a covering with translates of a
convex disk $K$ is by using the largest negatively-homothetic copy of $K$
instead of a circle. Then the problems of close packing and loose covering
become equivalent. He proved that, measuring closeness and looseness this
way, the triangle is the worst for both close packing and loose covering,
with closeness $2$ and looseness $3$. This result corresponds to the
theorems of F\'ary concerning the ``worst case'' for densest packing and thinnest
covering with translates of a convex disk, where again the triangle is
the worst one in both cases. Linhart also conjectured that the worst
case among centrally symmetric convex disks is the affine-regular octagon.
Specifically, he conjectured that every centrally symmetric convex disk
can pack the plane by translates with closeness at least $3+2\sqrt2$, and
can cover the plane by translates with looseness at least $4+2\sqrt2$.

{\sc{Zong}}~\cite{Zong08} considered the problem of closest packing,
with closeness measured the Linhart way, though phrased in a slightly
different manner, and he confirmed Linhart's conjecture about the
affine regular octagon. The relation between Linhart's approach to closeness
and Zong's so-called {\it simultaneous packing and covering constant} is
as follows: For a given convex disk $K$, let $c(K)$ denote the
minimum closeness (in Linhart's sense) of a packing with
translates of $K$, and let $\gamma(K)$ denote Zong's minimum
homothety coefficient for a transition from a packing to a covering
with translates of $K$. Then $\gamma(K)=1+c(K)^{-1}$.

Similarly to his simultaneous packing and covering
constant {\sc Zong} \cite{Zong03} introduced
the simultaneous lattice packing and covering constant $\gamma^*(K)$ and
he proved that $\gamma^*(K)\le7/4$ for every three-dimensional convex body
$K$.

Confirming a conjecture of {\sc L.~Fejes T\'oth} \cite{FTL76},
{\sc{B\"or\"oczky}} \cite{Boroczky86} proved that the closest packing
with unit balls in space is unique and is
obtained by placing the centers of the balls at the vertices and at the
centers of all cubes of a cubic lattice of edge-length $4/\sqrt {3}$.  Since
for balls the problems of closest packing and loosest covering are equivalent,
B\"or\"oczky's solution settles both. {\sc B\"or\"oczky} \cite{Boroczky01}
defined {\it{edge-closeness}} of a packing of congruent balls as the supremum
of the distance of a line of an edge of some Dirichlet cell and the corresponding
center of ball. He showed that the minimum edge-closeness of a packing of unit
balls is $\sqrt{3/2}$.  Again, the unique optimal arrangement is the
body-centered cubic lattice.

Using techniques of {\sc{Delone}} and {\sc{Ry\v{s}kov}} \cite{DeloneRyskov}
and {\sc{Ry\v{s}kov}} and
{\sc{Baranovski\u\i}} \cite{RyskovBaranovskii75,RyskovBaranovskii76},
{\sc{Hor\-v\'ath}} \cite{Horvath80,Horvath86} solved the problem of the
closest lattice packing with unit balls in dimension $4$ and $5$.
{\sc{Sch\"{u}rmann}} and {\sc{Vallentin}} \cite{SchurmannVallentin}
designed an algorithm for approximating the loosest lattice sphere covering
with arbitrary accuracy. In $6$-dimensional space, the algorithm produced the
best known lattice for loose sphere covering.

\section{Arranging regular tetrahedra}

Since no integer multiple of the dihedral angle $\varphi=\arccos(1/3)=1.23\ldots$
at the edges of the regular tetrahedron $T$ equals $2\pi$ ($5\varphi = 6.15\ldots$
is just slightly smaller than $2\pi$), we know that space cannot be tiled with
congruent copies of $T$, hence $\delta(T)< 1$. Then, as {\sc{Hilbert}}
\cite{Hilbert} asked, how densely can space be packed with congruent regular
tetrahedra? The question is also of interest in areas other than mathematics,
{\it e.g.}, physics (compacting loose particles), chemistry (material
design), etc. The past few years brought an exciting development: A series of
articles appeared, each providing a surprisingly dense---denser than any
previously known---packing. {\sc Conway} and {\sc{Torquato}}
\cite{ConwayTorquato} presented a packing with density $0.717455\ldots$,
which is almost twice the lattice packing density of the tetrahedron.
After improvements by {\sc{Chen}} \cite{ChenER},
{\sc{Haji-Akbari, Engel, Keys}}, {\sc{Zheng}}, {\sc{Petschek}},
{\sc{Palffy-Muhoray}} and {\sc{Glotzer}}
\cite{Haji-Akbari+}, {\sc{Kallus}}, {\sc{Elser}} and
{\sc{Gravel}} \cite{KallusElserGravel}, and {\sc{Torquato}} and
{\sc{Jiao}} \cite{TorquatoJiao09c}, a packing of the currently highest known
density, namely $4000/4671= 0.856347\ldots$ was  obtained by {\sc{Chen}},
{\sc{Engel}} and {\sc{Glotzer}} \cite{ChenEngelGlotzer}.

While it was known for a long time that the value of $\delta(T)$ must be
smaller than $1$, no explicit non-trivial ({\it i.e}., strictly below $1$)
upper bound for the packing density of the regular tetrahedron was presented
until {\sc Gravel}, {\sc Elser} and {\sc Kallus} \cite{GravelElserKallus} gave
such a bound, about $1-10^{-24}$. They also gave a similar upper bound for
the packing density of the regular octahedron. The gap between the density
of the best known packing and the upper bound remains quite wide, and it may be
very hard to narrow it down substantially. On the other hand, there is hope
for the determination of the translational packing density of tetrahedra for
which {\sc{Zong}} \cite{Zong19} suggested a computer approach.

The papers by {\sc{Lagarias}} and {\sc{Zong}} \cite{LagariasZong} and
{\sc{Ziegler}} \cite{Ziegler10} survey the history of packing regular tetrahedra.

The thinnest known covering by regular tetrahedra constructed by {\sc Conway} and
{\sc{Torquato}} \cite{ConwayTorquato} has density $9/8$.
{\sc{Fiduccia, Forcade}} and {\sc{Zito}} \cite{FiducciaForcadeZito}
and independently {\sc{Dougherty}} and {\sc{Faber}} \cite{DoughertyFaber} found
a body $T_{84}$ that admits a lattice tiling of $E^3$ and is inscribed in a
tetrahedron $T$ of volume ${\rm{vol}}(T)=\frac{125}{63}{\rm{vol}}(T_{84})$. It
follows that $\vartheta_L(T)\le\frac{125}{63}$. {\sc{Forcade}} and {\sc{Lamoreaux}}
\cite{ForcadeLamoreaux} conjectured that $\vartheta_L(T)=\frac{125}{63}$ and
support the conjecture by proving that the lattice corresponding to the tiling by
copies of $T_{84}$ represents a local minimum of the density. There is no
non-trivial lower bound for the covering density of the regular tetrahedron, but
recently {\sc{Xue}} and {\sc{Zong}} \cite{XueZong} established such a bound for
lattice arrangements: The lattice covering density of a simplex in $E^n$ satisfies
$$\vartheta_L(S)\ge1+\frac{1}{2^{3n+7}}.$$

\section{Packing cylinders}

The first known non-tiling 3-dimensional convex solid to have its packing density
determined was the infinite circular cylinder $C$, that is, the Minkowski sum
$B+L$ where $B$ is a circle and $L$ is a line, see {\sc{A.~Bezdek}} and
{\sc{W.~Kuperberg}} \cite{BezdekAKuperberg90}.  As expected,
$\delta(C)=\pi/\sqrt{12}$, the maximum density being attained
when all cylinders are parallel.

Concerning packings of unit cylinders in which no two are parallel,
{\sc{A.~Bezdek}} and {\sc{W.~Kuperberg}} \cite{BezdekAKuperberg91a} showed that
for every such packing, the complement of the packing contains a ball of radius
$r>\rho=\frac{2}{\sqrt3}-1$. In other words, the closeness (see XX) of such a
packing is greater than $1/\rho$. They also proved that every point $p$ lying in
the complement of the packing is within a distance of at most $\sqrt{2/3}$ from
the center of such a ball. Moreover, every ball of radius $\rho$ not
intersecting any of the cylinders can be moved continuously from its original
position to the position of any other such ball, while avoiding every cylinder
during the motion. It appeared that if no two of the cylinders in a packing
are parallel, then the density of the packing should be rather low, perhaps
even zero. However, {\sc K.~Kuperberg} \cite{KuperbergK} constructed such a
packing with positive density. {\sc Graf} and {\sc{Paukowitsch}}
\cite{GrafPaukowitsch} improved the construction, reaching density $5/12$.
By a further improvement {\sc{Ismailescu}} and {\sc{Laskawiec}}
\cite{IsmailescuLaskawiec} reached density $1/2$. Moreover, they constructed
packings of congruent cylinders with no two cylinders parallel to each other
whose local density in a ball of sufficiently large radius is arbitrarily
close to $\pi/\sqrt{12}$.

\section{Obstructing light}

H.~Hornich posed the question of how many material unit balls (meaning closed
balls with mutually disjoint interiors) are needed to radially shield one such
ball, in the sense that every ray emanating from the center of the shielded ball
must meet a shielding one. The set of shielding balls is called a {\it cloud}.
Let $H(r)$ denote the smallest number of unit balls in a cloud for a ball of
radius $r$.  As a simple consequence of (V.1.2) (page XX114), {\sc{L.~Fejes
T\'oth}} \cite{FTL59e} showed the inequality
$$
H(r)>\frac{12\alpha}{6\alpha-\pi},\quad
\frac{\pi}{6}<\alpha=\arctan{\frac{1+r}{\sqrt{6r+3r^2}}}<\frac{\pi}{2},
$$
which yields the lower bound $H(1)\ge19$. The bound was subsequently improved
to $H(1)\ge24$ by {\sc Heppes} \cite{Heppes67b} and then raised again by
{\sc Cs\'{o}ka} \cite{Csoka} to $H(1)\ge30$.  By a suitable construction,
{\sc Danzer} \cite{Danzer60} showed that $H(1)\le42$.

{\sc L.~Fejes T\'oth} \cite{FTL59e} showed that a point, and therefore also a
sufficiently small ball, can be shielded by six unit balls.
{\sc{Gr\"unbaum}} \cite{Grunbaum60} proved that five balls do not suffice.
This can be interpreted as the equality $H(0)=6$.  Besides $r=0$, the value of
$H(r)$ is not known for any $r$.  However, the above inequality yields
$H(r)>8\pi r^2/\sqrt{27}\,$. For large values of $r$ the estimate is
asymptotically exact.

The notion of a cloud for a ball can be extended in various ways. We can
consider clouds for balls of infinite radius, that is half-spaces. A
{\it{cloud for a half-space}} is a packing in the complement
of the half-space such that every line perpendicular to the bounding plane
intersects a member of the packing. Of course, a cloud for a half-space
contains infinitely many balls, so in this case, we are looking for the
{\it{width of the cloud}}, that is for the minimal width of a slab
containing the cloud. A packing disjoint from a set $S$ that intersects
all rays issuing from the boundary of $S$ in the complement of $S$ is a
{\it{dark cloud}} for $S$. If each such ray intersects the interior of a
member of the cloud then the cloud is called {\it{deep}}. Finally, if the
corresponding rays intersect at least $k$ members of the cloud, then we
have a {\it{$k$-fold cloud}}.

\medskip
\centerline {\immediate\pdfximage width5cm
{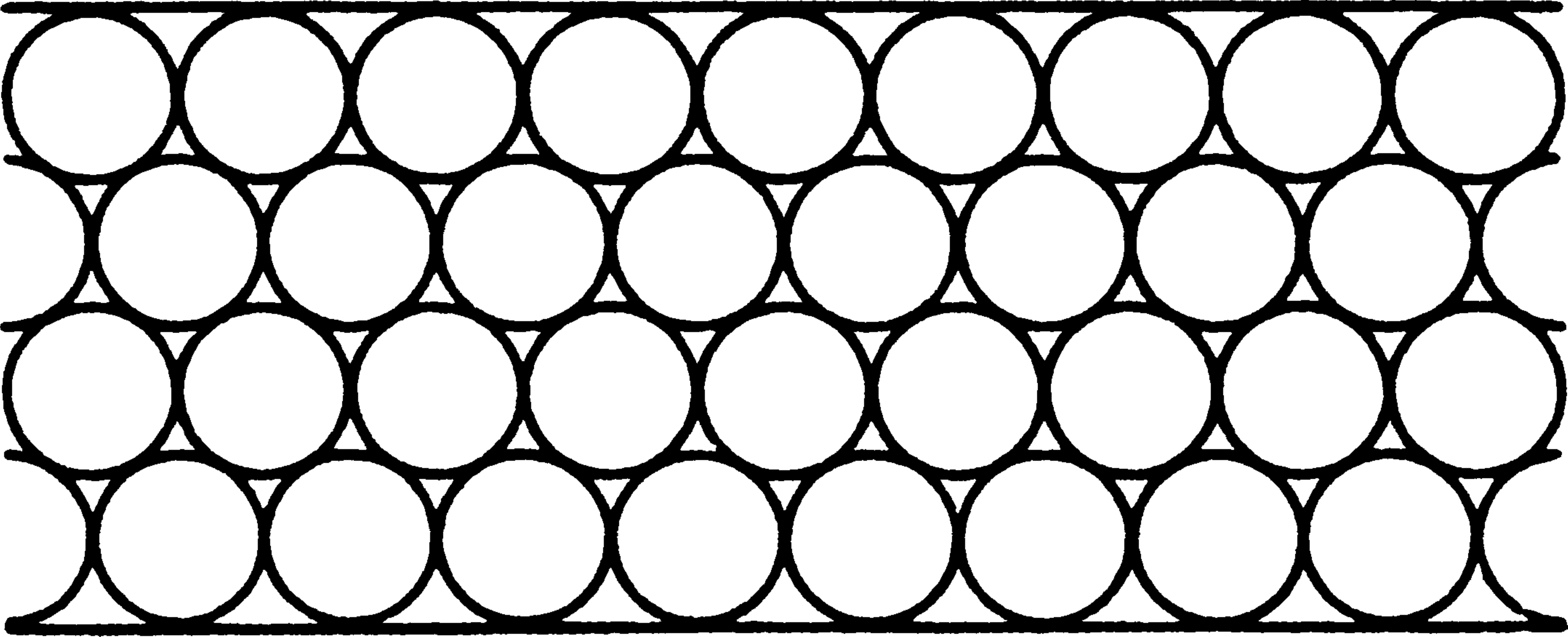}\pdfrefximage \pdflastximage}
\smallskip{\centerline{Figure~8}}
\medskip

Consider $k$ consecutive rows of the densest lattice packing of unit
circles. These circles are contained in a strip of width $(k-1)\sqrt3-2$
and form a $k$-fold cloud for the half-plane parallel to the rows.
{\sc{Heppes}}~\cite{Heppes61a} proved that this is the narrowest
$k$-fold cloud of unit circles shielding a half-plane. An alternative
proof of this statement is due to {\sc{Haj\'os}}~\cite{Hajos64}.
Figure~8 shows the narrowest $4$-fold cloud of circles.

It is not difficult to show ({\sc L.~Fejes T\'oth} \cite{FTL59e}) that the
width of any cloud of unit balls shielding a half-space in $E^3$ must be at least
$2+\sqrt2$. Equality holds only when the cloud consists of two square-pattern
layers in contact with each other, so that each ball in one layer touches
exactly four balls in the other layer. If we stack $k$ horizontal clouds of
this kind, we obtain a {\it $k$-fold cloud} of width $(2k-1)\sqrt2 + 2$.
{\sc Heppes} \cite{Heppes61a} improved this for every $k>1$ by
constructing a $k$-fold cloud of unit balls of width
$\left(k+\left\lfloor\frac{k-1}{3}\right\rfloor\right)\sqrt3 + 2$.
{\sc{Jucovi\v{c}}}~\cite{Jucovic66} gave bounds for the number of circles
in a $k$-fold cloud for a point.

For a convex body $K$, let $C(K)$ and $C_T(K)$ denote the minimum cardinality
of dark clouds for $K$ consisting of congruent copies of $K$ and of
translates of $K$, respectively.

The papers by {\sc{B\"or\"oczky}} and {\sc{Soltan}} \cite{BoroczkySoltan},
{\sc{Zong}} \cite{Zong97b, Zong99a}, {\sc{Talata}} \cite{Talata00a} contain
bounds for $C_T(K)$ for convex bodies $K$ in $E^n$. {\sc{Zong}} \cite{Zong97b}
proved that $C_T(K)\ge2n$ with equality only for parallelepipeds. Finding upper
bounds for $C_T(K)$ is related to the following problem:
Given $\varepsilon>0$, find the minimum number
$s(n,\varepsilon)$ such that for some packing of unit balls in $n$
dimensions, every line segment of length greater than $s(n,\varepsilon)$ must
be closer than $\varepsilon$ to one of the centers of the balls. The bounds
for $s(n,\varepsilon)$ given by {\sc Henk} and {\sc Zong} \cite{HenkZong00} were
improved by {\sc B\"or\"oczky, Jr.}~and {\sc Tardos} \cite{BoroczkyJrTardos},
whose result implies that
$$C_T(K)\le3^{{n^2}+o({n^2})}$$
for every convex body $K\in E^n$ and
$$C_T(K)\le2^{{n^2}+o({n^2})}$$
for every centrally symmetric convex body $K\in E^n$.

{\sc{Szab\'o}} and {\sc{Ujv\'ary-Menyh\'art}} \cite{SzaboUjvary-Menyhart}
proved that the minimum cardinality of a deep cloud for a convex disk is
at most 9 with equality only for the circle.

{\sc Heppes} \cite{Heppes60b} observed that a lattice packing of balls is
always penetrable by lines in three linearly independent directions, implying
that no lattice packing of balls can form a dark cloud for a half-space. This
shows that the existence of a dark cloud of congruent balls is not so trivial.
Nevertheless, {\sc B\"or\"oczky} \cite{Boroczky67} succeeded in constructing
such a cloud of a relatively small width.  Take four consecutive
hexagonal-pattern layers of balls, $s_0, s_1, s_2, s_3$ from the densest
lattice packing of balls.  B\"or\"oczky proved that these layers, adjoined by
the mirror images $s_{-1}, s_{-2}, s_{-3}$ of $s_1, s_2, s_3$ reflected in the
middle plane of $s_0$, form a dark cloud, and he asserted without proof that
the same is true for $s_{-2}, s_{-1}, s_0, s_1, s_2$ already.

The aforementioned theorem of Heppes was strengthened by {\sc{Hortob\'agyi}}
\cite{Hortobagyi71}, who proved that every lattice packing of
unit balls can be penetrated in three linearly independent directions by a
cylindrical beam of light of radius $\displaystyle {\frac{3\sqrt2}{4}}-1 =
0.06065\ldots$.  As the example of the densest lattice packing of balls shows,
this constant cannot be replaced with a larger one.  This last statement is
quite non-trivial, but it follows from Theorem II of {\sc Bambah} and {\sc
Woods} \cite{BambahWoods94}, see also {\sc Makai} and {\sc Martini}
\cite{MakaiMartini}, where a gap in the proof of this theorem was
filled.

The observation of Heppes inspired further research. {\sc{Hausel}} \cite{Hausel},
{\sc{Henk}} and {\sc{Zong}} \cite{HenkZong00}, {\sc{Henk}}, {\sc{Ziegler}} and
{\sc{Zong}} \cite{HenkZieglerZong} and {\sc{Henk}} \cite{Henk} estimated
the greatest number $k(n)$ with the property that for every lattice packing
of the $n$-dimensional ball there exists a ``free  $k$-dimensional plane'',
that is, a $k$-dimensional plane contained in the complement of the packing.
{\sc{Hausel}} \cite{Hausel} proved that $k(n)\le n-c\sqrt n$, for some
constant $c>0$. The best lower bound up to date, given by {\sc{Henk}}
\cite{Henk}, is $n/{\rm log}_2n\le k(n)$, for $n$ sufficiently large.
{\sc{Horv\'{a}th}} and {\sc{Ry\v{s}kov}} \cite{HorvathRyskov75a,HorvathRyskov75b}
estimated the maximum radius of a cylinder around a line that can be
inserted in the void of every lattice packing of the $n$-dimensional unit
ball, and conjectured that their result is sharp for $n=4$. Their conjecture was
refuted by {\sc{Makai}} and {\sc{Martini}} \cite{MakaiMartini}.

A related problem, posed by {\sc{G.~Fejes T\'{o}th}} \cite{FTG76b}, concerns
the thinnest lattice arrangement of balls that intersects every $k$-dimensional
plane. For $k=0$ the problem is about the thinnest lattice covering,
hence it is solved for $n\le5$. For $n=2$, $k=1$ the problem was solved by
{\sc L.~Fejes T\'{o}th} and {\sc Makai} \cite{FTLMakai} by proving that the
density of a lattice packing of circles intersecting every line is at least
$\sqrt3\pi/8$. For $k=n-1$ the problem turned out to be equivalent to finding
the densest lattice packing of balls. This is a consequence of the
following result of {\sc Makai} \cite{Makai78} also found independently
{\sc{Kannan}} and {\sc Lov\'{a}sz} \cite{KannanLovasz}: Let $\rho(K)$ denote the
infimum of the density of a lattice arrangement of a convex body $K$, such
that every hyperplane intersects one of the members of the arrangement, and
let $\widehat K$ denote the polar body of $\frac{1}{2}(K-K)$. Then
$$
\rho(K)\delta_L({\widehat K})={\rm vol}(K){\rm vol}\left(\frac{1}{4}{\widehat
K}\right).
$$
This solves the problem for balls in dimensions $n\le 8$ and $n=24$.  For
$0<k<n-1$ only the case $k=1$, $n=3$ has been solved: {\sc{Bambah}} and
{\sc{Woods}} \cite{BambahWoods94} showed that the thinnest lattice arrangement
of balls intersecting every line arises from the densest lattice packing
by enlarging the balls' radius by a factor of $3\sqrt2/\pi$. {\sc Kannan}
and {\sc Lov\'{a}sz} \cite{KannanLovasz} and {\sc{Gonz\'alez Merino}} and
{\sc{Schymura}} \cite{GonzalezMerinoSchymura} investigated the case
$0<k<n-1$ further.

\section{Avoiding obstacles}

If convex disks are packed in a parallel strip of the plane, we
call it a {\it layer of disks}. Let $w$ be the width of the strip and
$l$ be the length of a path that connects the two edges of the strip
without penetrating any of the disks.  ({\sc{L.~Fejes T\'oth}}
\cite{FTL66b}) defined the {\it permeability}
$p$ of the layer as
$$p=w/\inf l.$$
He showed that for the permeability $p$ of a layer of congruent
circles $p>\sqrt{27}/2\pi=0.82699\ldots$ holds. For layers consisting
of a large number of rows from the densest lattice packing of
circles, $p$ comes arbitrarily close to this lower bound. For
incongruent circles the value of $\inf p$ is not known. However,
by a rather complicated construction in the same article a layer
of incongruent circles with $p=0.82322\ldots<2\pi/\sqrt{27}$ was found.

The situation is quite different if one considers layers of squares
instead of circles. For layers consisting of congruent squares,
$\inf p=2/3$ holds. But the same holds even for layers of squares of
arbitrary sizes and orientations. Moreover, it was shown by
{\sc{L.~Fejes T\'oth}}~\cite{FTL68a} that to any layer of similar
copies of a parallelogram $P$ there is a layer consisting of
translates of a replica of $P$ with the same permeability.
It is an interesting question to consider which convex disks share
this property with the parallelogram and which behave like the
circle. {\sc{Florian}} and {\sc{Groemer}} \cite{FlorianGroemer}
showed that for every $m\ge39$ regular $m$-gons belong to the latter
group, that is there exists a layer of homothetic copies of them
whose permeability is smaller than the infimum of the permeability
of all layers of congruent regular $m$-gons.

The above mentioned results of L.~Fejes T\'oth were sharpened by
{\sc{Bollob\'as}} \cite{Bollobas68a} and {\sc{Florian}}
\cite{Florian78b}, Bollob\'as determining the infimum of the
permeability of a layer of given width of similar copies of a given
parallelogram and Florian determining the infimum of the
permeability of a layer of given width of unit circles.
{\sc{Hortob\'agyi}}~\cite{Hortobagyi76a} proved
that the permeability of a layer of translates of a disk of
constant width is at least $\sqrt{27}/2\pi$. {\sc{Florian}}
\cite{Florian79b} proved that the permeability of a layer
of translates of a regular hexagon, and also of a layer of
translates of a regular triangle is at least $3/4$. Subsequently,
{\sc{Florian}} \cite{Florian80a} observed that these results are
special cases of a general theorem. Namely, the infimum of the
permeability of layers by translates of a convex disk $K$
equals the infimum of the permeability of layers by translates
of the difference body $K-K$. A survey on permeability
can be found in {\sc{Florian}}~\cite{Florian80b}.

{\sc{L.~Fejes T\'oth}}~\cite{FTL78b} also raised the following, natural problem.
Given a packing of the plane with convex ``obstacles'' and two points not in
the interior of any of the obstacles, find or estimate the greatest possible
size of a necessary detour caused by the obstacles in traveling from one point
to the other. The first result in this direction was given by {\sc{Pach}}
~\cite{Pach77}, who proved that for any packing of square obstacles of sides at most
$1$, any pair of points at distance $d$ outside the squares can be connected
by an obstacle-avoiding path of length at most $\frac{3}{2}d+\sqrt{d}+1$.
{\sc{G.~Fejes T\'oth}}~\cite{FTG78} improved this bound to $(3d+1)/2$ and
proved the bound $(2\pi/\sqrt{27})(d-2)+\pi$ for unit circular obstacles.
Both of these bounds are sharp for infinitely many values of $d$.

The shortest path problem for balls was treated in {\sc{G.~Fejes T\'oth}}
\cite{FTG13}. It turned out that in $E^n$, even for a packing of balls with
arbitrary but bounded radii, where the obstacles' density might be 1, we need
not make a detour greater than $O(d/n)$ ($d\to\infty$) in order to connect
two points lying at distance $d$ outside the balls by a path avoiding the
balls. For a packing of congruent balls in $E^n$ the detour we have to make
approaches zero exponentially with the dimension.

Algorithmic aspects of problems of this type have been studied by
{\sc{Papadimitriou}} and {\sc{Yannakakis}}~\cite{PapadimitriouYannakakis89,PapadimitriouYannakakis91},
{\sc{Chan}} and {\sc{Lam}}~\cite{ChanLam}, and {\sc{A.~Bezdek}}~\cite{BezdekA99}.
A.~Bezdek also obtained a bound for the length of the shortest path in space with
cubical obstacles.

An interesting variation of the obstacle-avoiding path problem was considered
by {\sc{L.~Fejes T\'oth}}~\cite{FTL93}: From a point outside the obstacles, one
tries to escape in any direction to a distance $d$ away from the point.  He
proved that for any set of convex, open obstacles an escape path exists of
length at most $d^2+\frac{1}{2}\ln d + c$, and that some sets of obstacles
require the length of $d^2-\left(\frac{\pi}{3}+\frac{3}{\pi}\right)d +c$.

A dual problem arises when we consider a family of convex disks
covering the plane and we wish to travel only within the part of the plane
covered at least twice. {\sc{G.~Fejes T\'oth}}~\cite{FTG13} stated the following
conjecture. If the plane is covered by a family of unit circles, then for any
two points, each covered at least twice, there is a path contained in the
multiply covered region, connecting one point with the other, and of length at
most $d\sqrt2 + c\,$, where $d$ is the distance between the points and $c$ is
a constant. {\sc{Baggett}} and {\sc{A.~Bezdek}}~\cite{BaggettBezdekA} confirmed the
conjecture in the case when the circles form a lattice covering. {\sc{Rold\'an-Pensado}}
~\cite{Roldan-Pensado} showed that two points at distance $d$ apart lying in the
multiply covered part of the plane can be connected by a path that remains in the
part of the plane covered at least twice and whose length is at most
$(\pi/3+\sqrt3)d+c$ for some constant $c < 17$. {\sc{A.~Bezdek}} and {\sc{Yuan}}
\cite{BezdekAYuan} investigated a related problem, where they measured the length
of the path by the number of passings from one circle into another.

\section{Stability}

We say that a packing of (not necessarily congruent) circles is {\it{stable}}
if every circle is immobilized by the others, {\it i.e.}, if
in every circle the central angle $\lambda$ based on the largest
``free'' arc, that is containing no contact point, is smaller than $\pi$. The
{\it instability} of the packing is defined as $\Lambda=\sup\lambda$,
and its {\it stability} as $\pi-\Lambda$. We wish to find a thinnest
packing among those with prescribed stability. The following theorem
of {\sc{L.~Fejes T\'oth}}~\cite{FTL60a} is related to this problem:

Given a stable packing of the plane with circles whose radii have a
positive lower bound and a finite upper bound. If $d$ is the density
and $\Lambda$ is the instability of the packing, then
$$
d\ge\frac{\pi}{n\tan\frac{\Lambda}{2} +\tan\frac{2\pi-n\Lambda}{2}}\,,
\quad n=\lfloor2\pi/\Lambda\rfloor.
$$

This bound is reached in the case of congruent circles whose centers form the
vertices of one of the tilings $\{3,6\}$, $\{4,4\}$, $\{6,3\}$, $(4,8,8)$, or
$(3,12,12)$. The corresponding instabilities are $\Lambda=\pi/3$,
$2\pi/3$, $3\pi/4$ and $5\pi/6$. Surprisingly, {\sc{B\"or\"oczky}}
\cite{Boroczky64} succeeded in constructing stable packings of density $d=0$ consisting of
congruent circles (Figure~9). His example shows that the above inequality
is sharp even in the limiting case $\Lambda=\pi$. {\sc{Kahle}} \cite{Kahle} modified
B\"or\"oczky's construction to obtain thin stable packings of congruent circles in
rectangles and regular hexagons.

{\sc Dominy\'{a}k}~\cite{Dominyak} gave bounds for the instability of stable circle
packings in the spherical and hyperbolic plane. His bounds are sharp, in several cases
characterizing regular and Archimedean tilings.

\medskip
\centerline {\immediate\pdfximage width8cm
{abb_160.pdf}\pdfrefximage \pdflastximage}
\smallskip{\centerline{Figure~9}
\medskip

{\sc De Bruijn}~\cite{de Bruijn54,de Bruijn55} proved that the members of a packing of
starlike sets in $E^n$ can be moved apart arbitrarily far from each other, each
translated continuously without overlapping with any of the others (Figure~10), unless
each of the sets' star-center is unique, and they all coincide in the packing (Figure~11).
He also proved that in a finite packing of convex disks, for every direction, there is a
disk that can be translated in that direction arbitrarily far from its original position
without intersecting other members of the packing. {\sc{L.~Fejes T\'oth}} and {\sc{Heppes }}
\cite{FTLHeppes63} independently discovered the same results, and the result about starlike
sets was also noticed by {\sc{Dawson}} \cite{Dawson}. In contrast to the case of the plane,
in space there are finite packings of convex bodies such that no single member can move
rigidly without disturbing the others. The example given by L.~Fejes T\'oth and Heppes
consists of 12 tetrahedra packed around a rhombic dodecahedron. They conjectured that the
tetrahedra alone, without the central dodecahedron, have the stated property. This was
confirmed by {\sc{Snoeyink}} and {\sc{Stolfi}} \cite{SnoeyinkStolfi}. {\sc{Shephard}}
\cite{Shephard70} constructed a packing of twelve centrally symmetric polyhedra with
the same property.

\medskip
\centerline {\immediate\pdfximage height3cm
{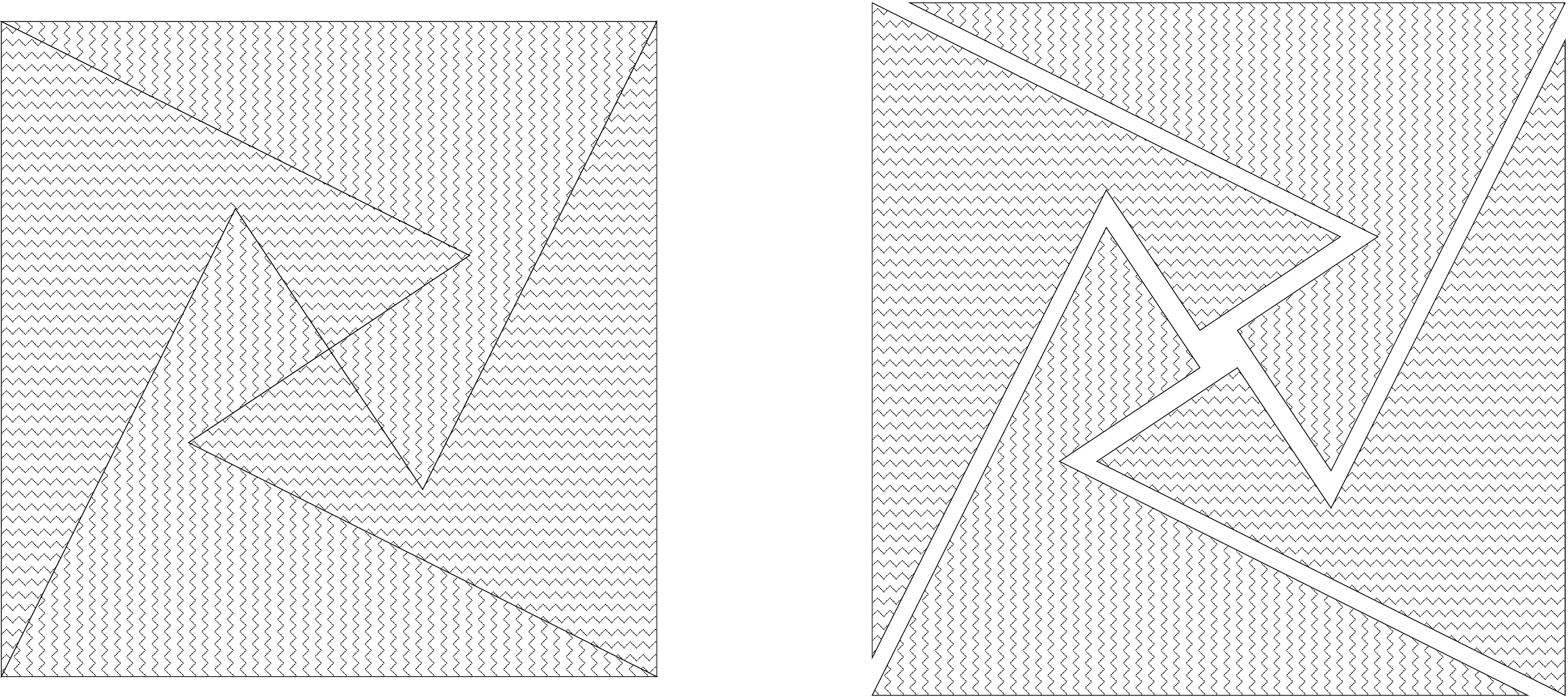}\pdfrefximage \pdflastximage\hskip1truecm
\immediate\pdfximage height3cm
{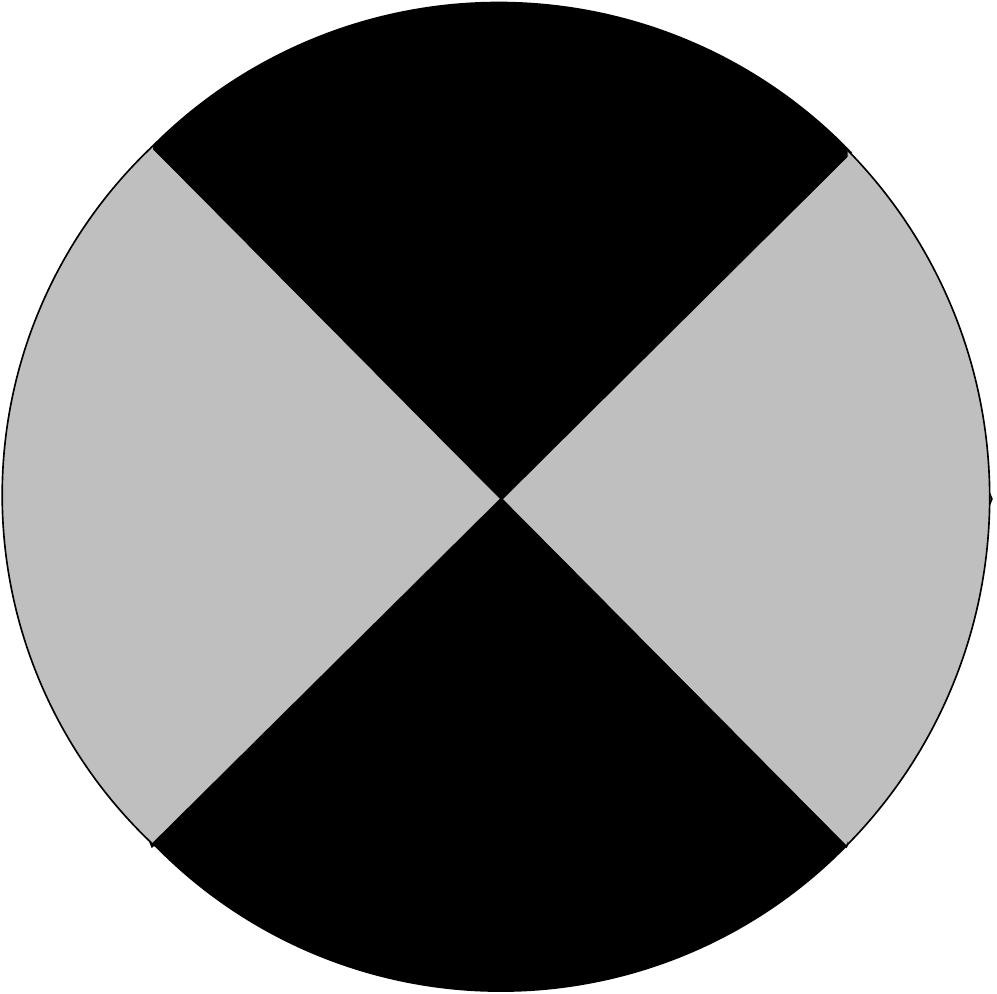}\pdfrefximage \pdflastximage}
\smallskip{\centerline{\hskip1.3truecmFigure~10 \hskip4truecm Figure~11}}
\medskip

By the theorem of De Bruijn, the members of every finite packing of convex bodies can be
moved arbitrarily far without disturbing the others by simultaneous translations.
{\sc{Natarajan}} \cite{Natarajan} conjectured that the members of such a packing can
even be {\it{separated by translations with two hands}} that is, a proper subset of the
packing exists that can be translated to infinity by applying a common translation to
them without disturbing the members in the complement. This conjecture turned out to
be false: {\sc{Snoeyink}} and {\sc{Stolfi}} \cite{SnoeyinkStolfi} gave a counterexample
of 6 bodies and showed that a packing of at most 5 convex bodies can indeed be separated
with two hands.

\section{Minkowskian arrangements}

A {\it Minkowskian arrangement} of similar copies of a centrally symmetric
convex disk was defined by {\sc{L.~Fejes T\'oth}} \cite{FTL65b} as an
arrangement in which no member contains the center of another one. He
proved that the density of such an arrangement of circles cannot exceed
$2\pi/\sqrt3$. A densest Minkowskian arrangement of circles consists of
congruent circles, and is obtained by replacing every circle in a densest
packing of congruent circles by a concentric one, with twice as large
radius. This result is an extensive generalization of inequality (III,2,1).
In fact, {\sc{L.~Fejes T\'oth}} \cite{FTL65b} proved that if finitely
many circles form a Minkowskian arrangement then the density of the
circles in their union cannot exceed $2\pi/\sqrt3$. In a subsequent paper
\cite{FTL67b} he gave an upper bound for the total area of such an
arrangement, which is sharp in many cases. {\sc{Moln\'ar}}~\cite{Molnar66a}
extended that result to the sphere and to the hyperbolic plane.

\medskip
\centerline {\immediate\pdfximage width5cm
{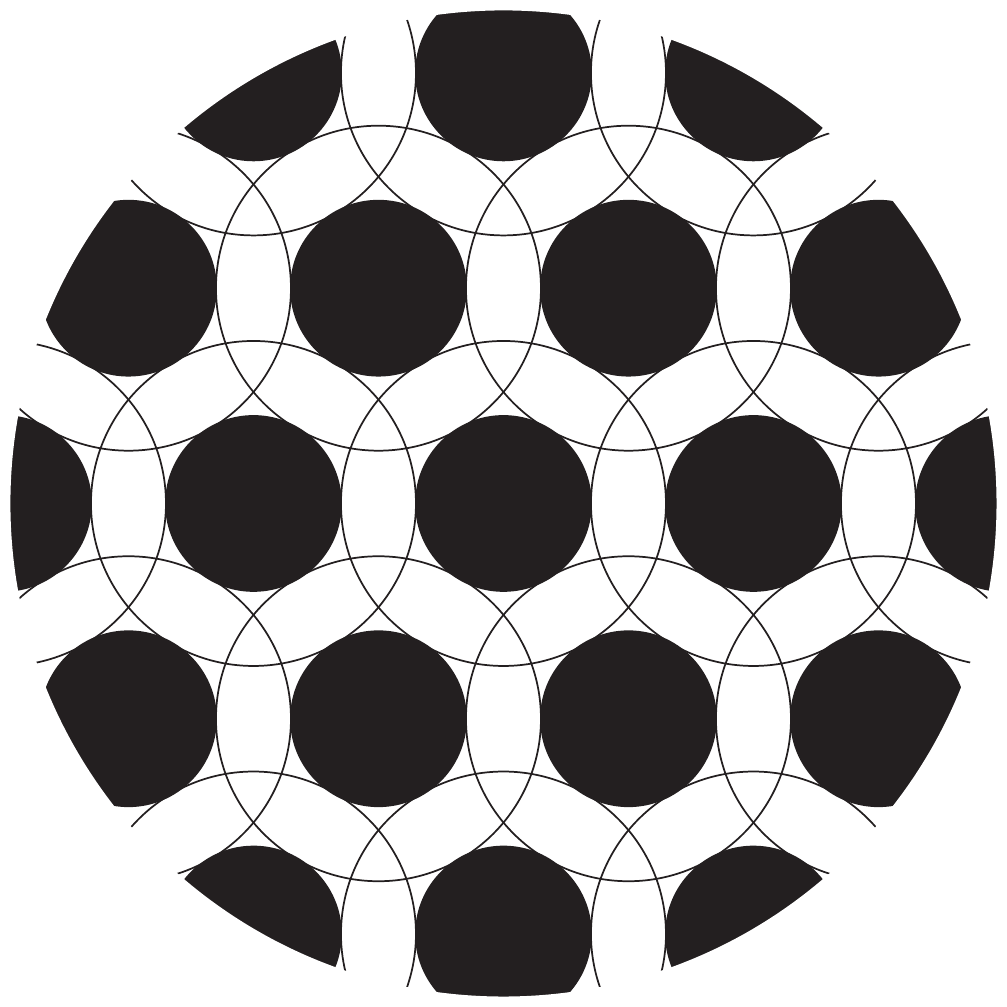}\pdfrefximage \pdflastximage}
\smallskip{\centerline{Figure~12}}
\medskip

The problem of the densest Minkowskian arrangement of circles was generalized by
{\sc{L.~Fejes T\'oth}}~\cite{FTL67a} as follows. Let $\mu$ be a positive number smaller than $1$.
We consider a set of circles $c_1, c_2\ldots$ of radii $r_1, r_2\ldots\,$. With each circle
$c_i$ we associate a concentric circle with radius $\mu r_i$, and we call it the {\it{kernel}} of
$c_i$. In a {\em{generalized Minkowskian arrangement of circles of order $\mu$}} none of the
circles is allowed to overlap the kernel of another. {\sc L.~Fejes T\'oth}~\cite{FTL67a}
conjectured that for $\mu\le\bar{\mu}=\sqrt3-1$ the densest arrangement consists of congruent
circles, and each of them touches six kernels (Figure~12). In this conjecture $\bar{\mu}$
denotes the greatest value of $\mu$ under which this particular arrangement is a covering.
{\sc Moln\'ar}~\cite{Molnar67b} and {\sc Florian}~\cite{Florian67} gave density
bounds under some condition for the homogeneity of the arrangement. {\sc{B\"or\"oczky}}
and {\sc{Szab\'o}}~\cite{BoroczkySzabo02} proved the conjecture in full generality.
{\sc{Kadlicsk\'o}} and {\sc{L\'angi}} \cite{KadlicskoLangi} gave a sharp bound for the total
area of the circles in a finite generalized Minkowskian arrangement of circles of order
$\mu\le\sqrt3-1$.

Minkowskian arrangements of circles on the sphere were treated by
{\sc{L.~Fejes T\'oth}} \cite{FTL99} who gave an upper bound for the total
area of $n$ spherical caps in a Minkowskian arrangement, sharp for
$n=3,\ 4,\ 6$ and $12$. The centers of the circles in the optimal
arrangement form an equilateral triangle inscribed in a great circle, a
regular tetrahedron, a regular octahedron and a regular icosahedron inscribed
in the sphere, respectively. Also, the bound is asymptotically sharp for
large $n$.

The density of a Minkowskian arrangement of homothetic copies of a centrally
symmetric convex body in $E^n$ is at most $2^n$. On the other hand, allowing
similar copies, there is no universal upper bound for the density. This was
noticed already by {\sc{L.~Fejes T\'oth}} \cite{FTL65b}, who also observed that
in order to achieve high density, the members of the arrangement must occur
in many different orientations: A Minkowskian arrangement of similar bodies in
$E^n$ with at most $m$ distinct orientations can have density at most $m2^n$.
{\sc{Bleicher}} and {\sc{Osborn}} \cite{BleicherOsborn} showed that there are
Minkowskian arrangements in $E^n$ even of congruent copies in at most $m$
distinct orientations with densities arbitrarily close to $m2^n$.

What is the maximum number $M(n)$ of pairwise intersecting homothetic copies
of a centrally symmetric convex body forming a Minkowskian arrangement in $E^n$?
The example of the cube shows that $3^n\le{M(n)}$. {\sc{F\"uredi}} and {\sc{Loeb}}
\cite{FurediLoeb}, who first considered this problem, proved that $M(n)\le 5^n$.
This upper bound was subsequently improved by {\sc{Nasz\'odi, Pach}} and
{\sc{Swanepoel}} \cite{NaszodiPachSwanepoel} to $O(3^n\ln{n})$, by
{\sc{Polyanskii}} \cite{Polyanskii17} to $3^{n+1}$, by {\sc{Nasz\'odi}} and
{\sc{Swanepoel}} \cite{NaszodiSwanepoel} to $2\cdot3^n$, and finally by
{\sc{F\"oldv\'ari}} \cite{Foldvari}, who proved the sharp bound $3^n$ with
equality only for parallelotopes.

\section{Saturated arrangements}

Minkowskian arrangements of circles are in a certain sense dual
to {\it saturated collections of circles}. Let $\mathcal{S}$ be a collection of
closed circular disks and let $r>0$ be the infimum of their radii. We say that
$\mathcal{S}$ is {\it saturated} if the part of the plane not covered by the circles
contains no circle of radius $r$. It was conjectured by
{\sc{L.~Fejes T\'oth}}~\cite{FTL67a} that the density of a saturated collection
of circles is always greater than or equal to $\pi/\sqrt{108}$. The conjecture
was confirmed by {\sc{Eggleston}}~\cite{Eggleston65} for families of mutually
non-overlapping circles, and proved in general by {\sc{Bambah}} and {\sc{Woods}}
\cite{BambahWoods68b}. Thus, a thinnest saturated arrangement of circles arises by
replacing the circles in a thinnest covering of the plane by concentric circles of half
size. This is a generalization of the inequalities (III,2,2) and (III,2,5). A corresponding
theorem for a packing of homothetic copies of a centrally symmetric convex disk in
place of circles was proved by {\sc{Bambah}} and {\sc{Woods}}~\cite{BambahWoods68a}.
{\sc{Dumir}} and {\sc{Khassa}} \cite{DumirKhassa73a} strengthened the above result
for circles, and in~\cite{DumirKhassa73b} for arbitrary centrally symmetric convex
disks as follows: No saturated arrangement of homothetic copies of a centrally
symmetric convex disk $K$ can cover a smaller portion of the plane than
$\vartheta(K)/4\,$. For the density of a saturated arrangement of homothetic copies of a
(not necessary symmetric) convex disk $K$ {\sc{Khassa}} \cite{Khassa75a} proved the
upper bound ${\rm{area}}(K)/t(K)$, where $t(K)$ denotes the area of the largest
triangle contained in $K$. For the density of saturated packings of balls in 3
dimensions {\sc{Khassa}} \cite{Khassa75b} established the upper bound $3/32$.

Recall from Chapter 10 that a packing with congruent copies of a set $K$ is
{\it $k$-saturated}, if deleting $k-1$ members of the packing never creates a void
large enough to pack in it $k$ copies of $K$. It is natural to ask for the infimum
$\Delta_k(K)$ of the densities of all $k$-saturated packings with replicas of $K$.
{\sc{G.~Fejes T\'oth, G.~Kuperberg}} and {\sc{W.~Kuperberg}}
\cite{FTGKuperbergKuperberg} proved the asymptotic bound
$\Delta_k(K)\ge\delta(K)-O(k^{-1/n})$ for every body $K$ in $E^n$. The determination
of these quantities is difficult even for $K=B^2$. The only result in this direction
is due to {\sc{Heppes}} \cite{Heppes01b}, who determined the infimum of the densities of
2-saturated lattice packings of circular disks, supporting the conjecture that
$\Delta_2(B^2)=\pi(3-\sqrt5)/\sqrt{27}=0.461873\ldots$

Another notion of higher order saturation was studied by {\sc{L.~Fejes T\'oth}} and
{\sc{Heppes}} \cite{FTLHeppes80} and {\sc{A.~Bezdek}} \cite{BezdekA90}. Here, we
formulate the concept only for the special case of packings of congruent circular disks.
A packing of disks of radius $r$ is {\it{saturated of order $k$}} if every disk of
radius $r$ intersects  at least $k$ members of it. Saturation of order 1 means just
saturation, and it is easily seen that the order of saturation of a packing of
congruent circles is at most 3. L.~Fejes T\'oth and Heppes proved that the density
of an order 3 saturated packing of congruent circles is at least $\pi/(2-\sqrt3)$,
and A.~Bezdek proved that the density of an order 2 saturated packing of congruent
circles is at least $\pi/(\sqrt{27})$. The thinnest order 3 saturated packing
arises by placing the centers in the vertices of a tiling $(3,3,4,3,4)$, and a
thinnest order 2 saturated packing consists of the face-incircles of the tiling
$\{3,6\}$.

\section{Compact packings}

A packing of the plane is said to be {\it compact} if each member $K$ of the
packing satisfies the following three conditions:


\begin{enumerate}
\item $K$ has a finite number of neighbors,
\item all neighbors of $K$ can be ordered cyclically so that
each of them touches its successor,
\item the union of the neighbors of $K$ contains a polygon
enclosing $K$.
\end{enumerate}
\smallskip

{\sc{L.~Fejes T\'oth}}~\cite{FTL84b} originated the study of compact packings by
proving that if a compact packing of the plane with circular disks has positive
homogeneity then its density is at least $\pi /\sqrt{12}$. Further, if a
compact packing of the plane with homothetic centrally symmetric convex disks
has positive homogeneity, then its density is at least $3/4$, where equality
occurs only for packings with affine regular hexagons. {\sc{A.~Bezdek, K.~Bezdek}}
and {\sc{B\"or\"oczky}}~\cite{BezdekABezdekKBoroczky} proved that if
a compact packing of the plane with positively homothetic copies of a convex
disk has positive homogeneity, then its density is at least $1/2$, and equality
occurs for various packings with homothetic triangles.

The only compact packing of congruent circular disks is the hexagonal lattice.
{\sc{Kennedy}}~\cite{{Kennedy06}} considered compact packings of circular disks of two
different radii, $1$ and $r<1$, and proved that there are only nine values of
$r$ for which such compact packings exist. He also described all packing
configurations in the nine cases. {\sc{Messerschmidt}} \cite{Messerschmidt20}
proved the upper bound 13617 for the number of pairs ($r,s)$ that allow a
compact packing by disks of radii 1, $r$ and $s$ ($r<s<1$). In fact, there are
much fewer such pairs: {\sc{Fernique, Hashemi}} and {\sc{Sizova}}
\cite{FerniqueHashemiSizova} enumerated all 164 compact packings consisting
of three different sizes of circular disks. {\sc{Messerschmidt}}
\cite{Messerschmidt21} proved that
for every $n$ there exist only finitely many tuples $(r_1,\ldots,r_n)$ with
$0< r_1<\ldots<r_n=1$ that can occur as the radii of the disks in any compact
packing of the plane with $n$ distinct sizes of disk.

It can be expected that compact packings of circles are the densest among
all packings with the given radii. {\sc{B\'edaride}} and {\sc{Fernique}}
\cite{BedarideFernique} conjectured that if disks with $n$ different radii
allow a saturated compact packing in which at least one disk of each size
appears, then the maximal density over all the packings by disks with
these $n$ radii is reached for a compact packing. The hypothesis of saturation
is necessary for $n\ge3$. The conjecture was proved for for some cases,
including all of the nine pairs of radii allowing compact
packings of disks with two different sizes, by {\sc{Heppes}} \cite{Heppes00,Heppes03b},
{\sc{Kennedy}} \cite{Kennedy04}, {\sc{B\'edaride}} and {\sc{Fernique19a}}
\cite{BedarideFernique} and {\sc{Fernique}} \cite{Fernique19a}, however {\sc{Fernique}}
and {\sc{Pchelina}} \cite{FerniquePchelina} showed that the density of one of the
compact circle packings with three different radii is not maximal.

{\sc{Florian}}~\cite{Florian85} considered compact packings with circular disks on the
sphere and in the hyperbolic plane.

Compact packings in higher dimensions were investigated by {\sc{K.~Bezdek}}
\cite{BezdekK87b} and {\sc{K.~Bezdek}} and {\sc{Connelly}} \cite{BezdekKConnelly91}.
A packing in $E^n$ is {\it{compact}} if each member $A$ of the packing is enclosed by
its neighbors in the sense that any curve, connecting a point of $A$ with a point
sufficiently far from $A$, intersects the closure of a neighbor of $A$. {\sc{K.~Bezdek}}
and {\sc{Connelly}} \cite{BezdekKConnelly91} proved that the density of a compact
packing in $E^n$ consisting of homothetic centrally symmetric convex bodies with bounded
homogeneity is at least $(n+ 1)/2n$, and there is a compact lattice packing of centrally
symmetric convex bodies where equality holds.

{\sc{Fernique}} \cite{Fernique21} proposed a different generalization of the concept of
compact packings. The contact graph of a compact packing of circles, i.e., the graph that
connects the centers of adjacent circles, is a triangulation. By analogy, Fernique calls
a packing of balls in $E^n$ compact if its contact graph is the 1-skeleton of a face-to-face
tiling by simplices. Since regular tetrahedra do not tile the space, there is no compact
packing of congruent balls. {\sc{Fernique}} \cite{Fernique21} determined the unique
compact packing with two different sizes of balls and in  \cite{Fernique19b} he
described all the four compact packings with balls of three different sizes.

\section{Totally separable packings}

A packing of convex bodies is {\it totally separable} if each pair of the
bodies can be separated by a hyperplane not intersecting the interior of any of the
bodies. Given a convex disk, what is the maximum density of a totally
separable packing with its congruent copies? {\sc{G.~Fejes
T\'oth}} and {\sc{L.~Fejes T\'oth}}~\cite{FTGFTL73b} proved that the density of an
arbitrary totally separable packing with congruent copies of a convex disk
cannot exceed the ratio between the area of the disk and the minimum area of
a quadrilateral containing the disk. If the disk is centrally symmetric,
then that ratio is actually the maximum density of such a packing. Namely,
the minimum  area of the quadrilateral containing a centrally symmetric disk
can always be attained by a parallelogram, and congruent parallelograms admit a
totally separable tiling. In particular, this yields that the density of a
completely separable packing with congruent circles cannot exceed $\pi/4$.
{\sc{A.~Bezdek}}~\cite{BezdekA83} proved this density bound under an assumption
weaker than total separability, requiring only the packing to satisfy the
following {\it local separability} condition: For every triple of circles
there is a line separating one of them from the other two. He also showed
that, in general, for non-circular disks local separability does not imply
the similar parallelogram bound. {\sc{K. Bezdek}} and {\sc{L\'angi}}
\cite{BezdekKLangi20b} proved an analogue of Oler's inequality for totally
separable packings of translates of a convex disk. {\sc{Vermes}}
\cite{Vermes96} investigated totally separable tilings and totally separable
packings of circles in the hyperbolic plane.

\bigskip
\centerline {\immediate\pdfximage width8truecm
{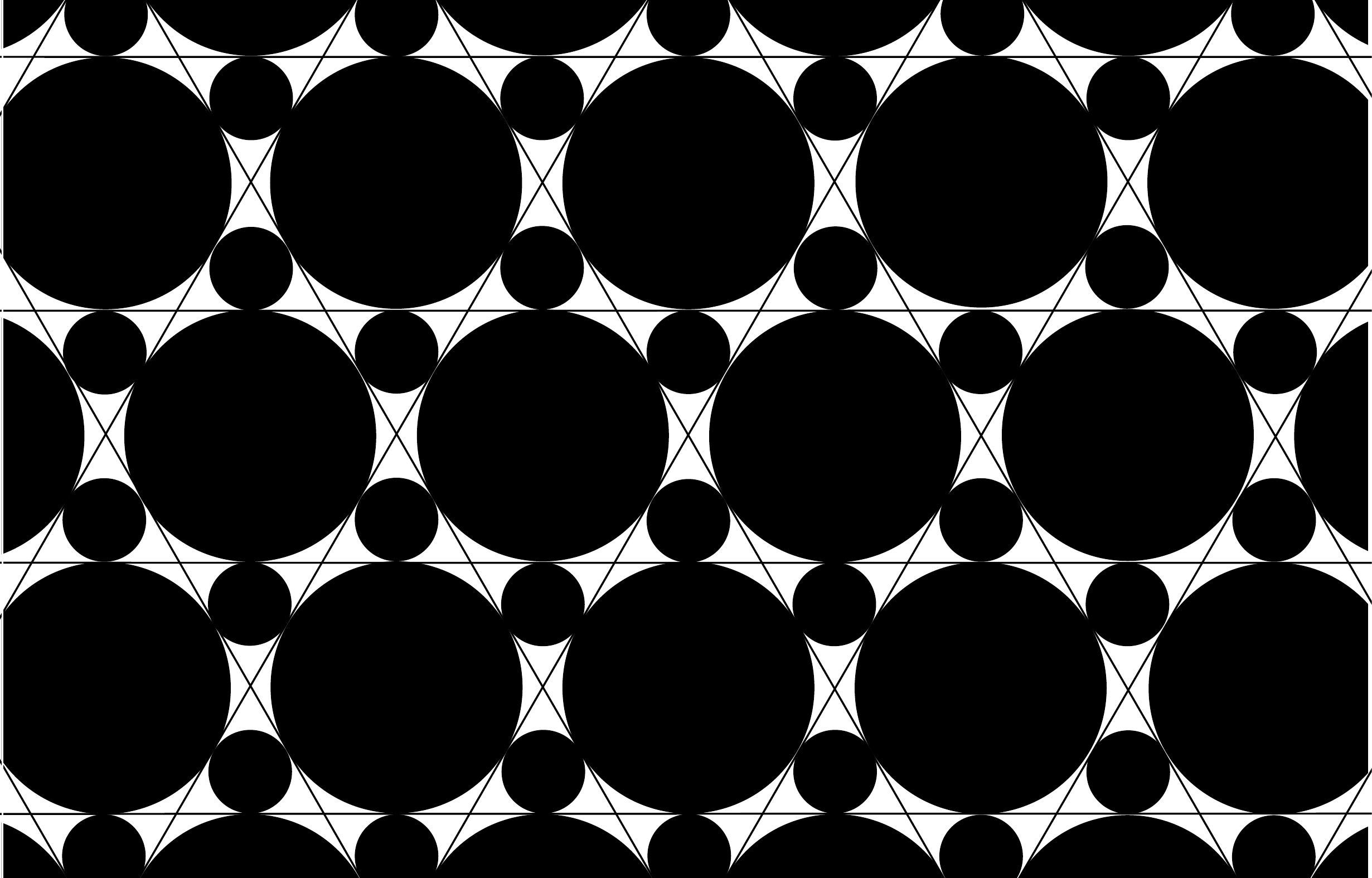}\pdfrefximage \pdflastximage}
\bigskip{\centerline{Figure~13}}
\bigskip

The interesting problem of maximum density of a totally separable
packing of (not necessarily congruent) circular disks remains
open. One can rephrase this problem in the following way: Start with a
configuration of lines partitioning the plane into bounded regions and
place a circle in each of the regions. What is the maximum density of a
circle packing so obtained? In this phrasing, the analogous question about
the minimum density of the covering by the circumcircles of the resulting
cells can be asked. It is conjectured that the configuration of lines
providing each of these extreme densities partitions the plane into the
Archimedean tiling $(3,6,3,6)$ (see Figure~13). {\sc{G.~Fejes T\'oth}}
\cite{FTG87} proved that neither of the extreme densities in question
are equal to $1$.

{\sc{K.~Bezdek}} and {\sc{L\'angi}} \cite{BezdekKLangi22} and {\sc{V\'as\'arhelyi}}
\cite{Vasarhelyi03} investigate totally separable packings and coverings
of circles on the sphere.

The problem of the densest totally separable packing of
$E^3$ with congruent balls was solved by {\sc Kert\'esz} \cite{Kertesz88}, who
proved that if a cube of volume $V$ contains $N$ unit balls forming a totally
separable packing, then $V \ge 8N$. Consequently, the cubic lattice packing is
the densest one among all totally separable packings of $E^3$ with congruent
balls, and the maximum density is $\pi/6$.

\section{Point-trapping lattices}

An arrangement of sets is {\it point-trapping} if every component of the
complement of the union of the sets is bounded. It is natural to ask: What is
the minimum density of a point-trapping lattice arrangement of any
$n$-dimensional convex body?  Confirming the {\it chessboard conjecture} of
{\sc L.~Fejes T\'{o}th} \cite{FTL75c}, {\sc B\"or\"oczky}, {\sc B\'ar\'any},
{\sc Makai} and {\sc Pach}  \cite{BoroczkyBaranyMakaiPach} proved that the
minimum is $1/2$, attained in the ``chessboard'' lattice arrangement of cubes.

The problem about the minimum density of a point-trapping lattice of $K$ can
be posed for any specific convex body $K$.  {\sc Bleicher} \cite{Bleicher75}
proved that the minimum density of a point-trapping lattice arrangement of
three-dimensional unit balls is equal to $128\pi/3{\sqrt{7142+1802\sqrt{17}}}=
1.1104\ldots$, and the extreme lattice is generated by three vectors, each of
length $\frac{1}{2}\sqrt{7+\sqrt{17}}=1.6676\ldots$ and each two forming an
angle of $\arccos\frac{\sqrt{17}-1}{8}=67.021\ldots^\circ$.

\section{Connected arrangements}

An arrangement of sets is said to be {\it connected} if the union of the sets
is connected. The problem of the minimum density $c(K)$ of a connected lattice
arrangement of an $n$-dimensional convex body $K$ has been explored by
{\sc Groemer} \cite{Groemer66b}, who proved the inequalities
$$
\frac{1}{n!}\le c(K)\le\frac{\pi^{n/2}}{2^n\Gamma(\frac{n}{2}+1)}.
$$
Both inequalities are sharp. The value $c(K) ={1/n!}$ is attained when
$K$ is a simplex or a cross-polytope, and the other extreme value of
$c(K)$ is attained when $K$ is a ball. Groemer characterized those
centrally symmetric bodies $K$ for which $c(K)=1/n!$. Extending Groemer's
investigation, {\sc L.~Fejes T\'oth} \cite{FTL73d} characterized all
$n$-dimensional convex bodies $K$ for which the inequality
$c(K)\ge1/n!$ turns into equality: They are the topological isomorphs of
the regular cross-polytope and their limiting polytopes.

\section{Points on the sphere}

In Chapter XX we mentioned results about problems of finding the optima of
sums of the form $\sum_{i\ne j}f(|x_i-x_j|)$ for a given number $N$ of points
$\{x_1,\ldots,x_N\}$ on the unit sphere. We continue to survey problems of
this type. We start with two problems of this kind, of interest from the
geometric point of view, and remarkably easy to solve for every value of $N\ge2$.

Distribute points $x_1,\ldots,x_N$ on $S^2$ so that:

\begin{enumerate}
\item  the sum $\sum(|x_i-x_j|^2$ is as large as possible;

\item  the sum $\sum\widehat{x_ix_j}$ is as large as possible, where
    $\widehat{x_ix_j}=2\arcsin\left(\frac{1}{2}|x_i-x_j|\,\right)$
     is the spherical distance between $x_i$ and $x_j\,$.
\end{enumerate}

Concerning the first problem we have $\sum(|x_i-x_j|^2\le N^2$
and equality occurs only if the vectors from the sphere's center to the points
$P_i$ are in equilibrium (see {\sc{L.~Fejes T\'oth}} \cite{FTL56b}).

For the solution of the second problem, the parity of $N$ plays a role. If
$N=2k$, then $\sum\widehat{x_ix_j}\le \pi k^2$, and equality occurs precisely
when the configuration of points is symmetric about the sphere's center. For
$N=2k+1$ we have $\sum\widehat{x_ix_j}\le\pi k(k+1)$. Equality occurs here
precisely when all points without an antipodal partner are distributed on a
great circle in such a way that two open semicircles determined by any point
contain the same number of points. The case $N=4$ was solved by
{\sc{Frostman}} \cite{Frostman}, the cases $N=5$ and $6$ by {\sc{L.~Fejes T\'oth}} \cite{FTL59d},
the cases $N=2k\ \ (k=1,2,\ldots)$ by {\sc{Sperling}} \cite{Sperling}, and the general
case by {\sc{Nielsen}} \cite{Nielsen} (see also {\sc{Larcher}} \cite{Larcher}).

{\sc{L.~Fejes T\'oth}} \cite{FTL59d} also considered the corresponding problem in
the elliptic plane, that is maximizing the sum of the non-obtuse angles formed
by $N$ lines. He conjectured that the maximum of the angle sum is attained when the lines are
evenly distributed between the three coordinate axis in which case the sum of the
angles is asymptotically $\frac{N^2\pi}{6}=N^2\cdot0.523\ldots$. He solved the problem
for $N\le6$ and proved the upper bound $\frac{k(N-1)}{5}\pi$ on the sum of the angles
of $N$ lines. {\sc{Fodor, V\'{\i}gh}} and {\sc{Zarn\'ocz}} \cite{FodorVighZarnocz16a}
gave an improvement of this bound that is asymptotically equal to
$\frac{3N^2\pi}{16}=N^2\cdot0.589\ldots$ as $N\to\infty$. {\sc{Bilyk}} and {\sc{Matzke}}
\cite{BilykMatzke} further improved the bound to
$(\frac{\pi}{4}-\frac{69}{300})N^2=N^2\cdot0.555\ldots$ and also gave a bound for
the corresponding problem in higher dimensions.

{\sc{Alexander}} \cite{Alexander72} considered another problem of this kind: Find the
configurations of $N$ points on the $n$-dimensional unit sphere for which the sum of all
distances between the points attains its maximum $S(N,n)$. Through an elegant
integral-averaging technique he obtained the bounds
$$\frac{2}{3}N^2-10\sqrt{N}\le S(N,2)\le\frac{2}{3}N^2-\frac{1}{2}\,.$$
{\sc Stolarsky} \cite{Stolarsky72}
extended Alexander's result for powers of the distances between the points and
generalized it to all dimensions. In \cite{Stolarsky73} {\sc Stolarsky} proved a
remarkable invariance theorem stating that the sum of all the distances between
the points of a given set of $N$ points on the $n$-dimensional sphere plus the
discrepancy of the set is a constant independent from the distribution of the points.
Using this result, Stolarsky gave a sharper bound for $S(N,n)$. If $c_0(n)$ denotes
the average distance from a variable point to a fixed one on the surface of the
sphere, then there is a constant $c_1(n)$ such that
$$
S(N,n)<c_0(n)N^2 - c_1(n)N^{1-1/n}.
$$
{\sc Beck} \cite{Beck} showed exactness of the constant $c_0$ by proving that
$$
c_0(n)N^2 - c_2(n)N^{1-1/n}<S(N,n)
$$
with a suitable constant $c_2(n)$. Alternative proofs of Stolarsky's invariance
theorem were given by {\sc{Brauchart}} and {\sc{Dick}} \cite{BrauchartDick} and
{\sc{Bilyk, Dai}} and {\sc{Matzke}} \cite{BilykDaiMatzke}. The latter authors
proved various generalizations of the invariance principle, and applied them to
problems of energy optimization.

In a similar vein, {\sc Witsenhausen} \cite{Witsenhausen} stated the following problem:
Under the constraint that the diameter of the set of points $x_1,x_2,\ldots,x_k$
in $n$-dimensional space must be smaller than $1$, find their configuration that maximizes
the sum $\sum\limits_{i,j} |x_i-x_j|^2$. Witsenhausen conjectured that the maximum,
denoted by $M(n,N)$, is attained when the points are distributed among the $n+1$
vertices of a regular simplex of edge-length $1$ and supported the conjecture by
proving the inequality $M(n,N)\le N^2n/(n+1)$. The conjecture was confirmed for
$n=2$ by {\sc Pillichshammer} \cite{Pillichshammer} and in arbitrary dimension by
{\sc{Benassi}} and {\sc{Malagoli}} \cite{BenassiMalagoli}.

\section{Arrangements of great circles}

{\sc{Maehara}} \cite{Maehara95} considered the following problem: Arrange $k$ great circles
so that the maximum spherical distance between a point of the sphere and the
nearest crossing point of the great circles is as small as possible. He solved
the problem for $k=3$ and $4$. The optimal arrangement in these cases occur
when each circle is divided into $2(k-1)$ equal arcs by the other $k-1$ great
circles.

{\sc{L.~Fejes T\'oth}} \cite{FTL59b} showed that the longest edge of a spherical tiling
formed by four or five great circles takes its minimum for the Archimedean tiling
(3,4,3,4) and (3,5,3,5), respectively. {\sc{Heppes}} \cite{Heppes58} studied a related
problem: He showed that the area of the face with the smallest area of a tiling formed
by four great circles takes its maximum for the cuboctahedron (3,4,3,4).

A conjecture of {\sc{L.~Fejes T\'oth}} \cite{FTL87} about arrangements of great
circles on a sphere states: If $k$ great circles are in general position (no
three of them have a common point), then the ratio of the greatest among
the areas of the regions into which the circles partition the sphere to the
smallest one tends to infinity as $k\to\infty$. In the same article a lower
bound is given: For sufficiently large $k$ the ratio is greater than $7.43$.

Motivated by this conjecture, {\sc{Ismailescu}} \cite{Ismailescu03} considered
an analogous problem on the plane and proved the following: Consider an
arrangement of $k$ lines such that no three are concurrent and all intersection
points lie inside a unit circular disk. Then among the $1+k(k-1)$ bounded
cells of the subdivision of the plane by the lines, there is one whose area
is at least $\frac{\pi}{4k}$. As a corollary, it follows that the ratio $q$
of the greatest area to the smallest area of cells is at least $(k+1)/8$. The
question whether the ratio $q$ tends to infinity as $k\to\infty$ if we do
not restrict the points of intersections to lie in a circle, remains open.


\small{
\bibliography{pack}}}
\begin{thebibliography}{2000}
\font \sc = cmcsc8
\parindent 0pt

\bibitem[]{}


\vskip2pt\hskip-1.5truecm{\sc{Alexander, R.}}

\bibitem[1972]{Alexander72} On the sum of distances between $n$ points on a sphere.
Acta Math. Acad. Sci. Hungar. {\bf{23}} (1972) 443--448. MR0312395, 
DOI~10.1007/BF01896964

\vskip2pt\hskip-1.5truecm{\sc{Baggett, D. R.}} and {\sc{Bezdek, A.}}

\bibitem[2003]{BaggettBezdekA} On a shortest path problem of G. Fejes T\'oth.
Discrete geometry, 19--26, Monogr. Textbooks Pure Appl. Math., 253, Dekker, New York, 2003. MR2034704 

\vskip2pt\hskip-1.5truecm{\sc{Bambah, R. P.}} and {\sc{Woods, A. C.}}

\bibitem[1968a]{BambahWoods68a} On the minimal density of maximal packings of the plane by convex bodies.
Acta Math. Acad. Sci. Hungar.  {\bf{19}}  (1968) 103--116. MR0223978, 
DOI~10.1007/BF01894686

\bibitem[1968b]{BambahWoods68b} On minimal density of plane coverings by circles.
Acta Math. Acad. Sci. Hungar. {\bf{19}} (1968) 337--343. MR0235462, 
DOI~10.1007/BF01894513

\bibitem[1994]{BambahWoods94} On a problem of G. Fejes T\'oth.
K. G. Ramanathan memorial issue. Proc. Indian
Acad. Sci. Math. Sci. {\bf{104}} (1994) no. 1, 137--156. MR1280062, 
DOI~10.1007/BF02830878

\vskip2pt\hskip-1.5truecm{\sc{Beck, J.}}

\bibitem[1984]{Beck} Sums of distances between points on a sphere ---
an application of the theory of irregularities of distribution to discrete geometry.
Mathematika {\bf{31}} (1984) no. 1, 33--41. MR0762175,
DOI~10.1112/S0025579300010639

\vskip2pt\hskip-1.5truecm{\sc{B\'edaride, N.}} and {\sc{Fernique, T.}}

\bibitem[2020]{BedarideFernique} Density of binary disc packings: the
9 compact packings.
To appear in \emph{Discrete Comput.\ Geom.}
arXiv:2002.07168v1 [cs.DM] 17 Feb 2020

\vskip2pt\hskip-1.5truecm{\sc{Benassi, C.}} and {\sc{Malagoli, F.}}

\bibitem[2008]{BenassiMalagoli} The sum of squared distances under a diameter constraint, in arbitrary dimension.
Arch. Math. (Basel) {\bf{90}} (2008) no. 5, 471--480. MR2414250,
DOI~10.1007/s00013-008-2509-z

\vskip2pt\hskip-1.5truecm{\sc{Betke, U.}} and {\sc{Henk, M.}}

\bibitem[1998]{BetkeHenk98} Finite packings of spheres.
Discrete Comput. Geom. {\bf{19}} (1998) no. 2, 197--227. MR1600046,
DOI10.1007/PL00009341

\vskip2pt\hskip-1.5truecm{\sc{Betke, U.}}; {\sc{Henk, M.}} and {\sc{Wills, J. M.}}

\bibitem[1994]{BetkeHenkWills} Finite and infinite packings.
J. Reine Angew. Math. {\bf{453}} (1994) 165--191. MR1285784,
DOI~0.1515/crll.1994.453.165

\vskip2pt\hskip-1.5truecm{\sc{Bezdek, A.}}

\bibitem[1980]{BezdekA80} Remark on the closest packing of convex discs.
Studia Sci. Math. Hungar. {\bf{15}} (1980) no. 1-3, 283--285. MR0681449 

\bibitem[1983]{BezdekA83} Locally separable circle packings.
Studia Sci. Math. Hungar. {\bf{18}} (1983) no. 2-4, 371--375. MR0787941 

\bibitem[1990]{BezdekA90} Double-saturated packing of unit disks.
Period. Math. Hungar. {\bf{21}} (1990) no. 3, 189--203. MR1105706, 
DOI~10.1007/BF02651088

\bibitem[1999]{BezdekA99} On optimal route planning evading cubes in the three space.
Beitr\"age Algebra Geom. {\bf{40}} (1999) no. 1, 79--87. MR1678579 

\vskip2pt\hskip-1.5truecm{\sc{Bezdek, A.}}; {\sc{Bezdek, K.}} and {\sc{B\"or\"oczky, K.}}

\bibitem[1986]{BezdekABezdekKBoroczky} On compact packings.
Studia Sci. Math. Hungar. {\bf{21}} (1986) no. 3-4, 343--346. MR0919378 

\vskip2pt\hskip-1.5truecm{\sc{Bezdek, A.}} and {\sc{Kuperberg, W.}}

\bibitem[1990]{BezdekAKuperberg90} Maximum density space packing with congruent circular cylinders of infinite length.
Mathematika {\bf{37}} (1990) no. 1, 74--80. MR1067888, 
DOI~10.1112/S0025579300012808

\bibitem[1991a]{BezdekAKuperberg91a} Placing and moving spheres in the gaps of a cylinder packing.
Elem. Math. {\bf{46}} (1991) no. 2, 47--51. MR1097930 

\bibitem[1997]{BezdekAKuperberg97} Circle covering with a margin.
3rd Geometry Festival: an International Conference on Packings, Coverings and Tilings (Budapest, 1996).
Period. Math. Hungar. {\bf{34}} (1997) no. 1-2, 3--16. MR1608299,
DOI~10.1023/A:1004275121

\vskip2pt\hskip-1.5truecm{\sc{Bezdek, A.}} and {\sc{Yuan, L.}}

\bibitem[2020]{BezdekAYuan} Short chains in circle and in square coverings.
Beitr. Algebra Geom. {\bf{61}} (2020) no. 1, 129--137. MR4058108,
DOI~10.1007/s13366-019-00462-x

\vskip2pt\hskip-1.5truecm{\sc{Bezdek, K.}}

\bibitem[1987a]{BezdekK87a} Densest packing of small number of congruent spheres in polyhedra.
Ann. Univ. Sci. Budapest. E\"otv\"os Sect. Math. {\bf{30}} (1987) 177--194. MR0927822

\bibitem[1987b]{BezdekK87b} Compact packings in the Euclidean space.
Beitr\"age Algebra Geom. {\bf{25}} (1987) 79--84. MR0899386 

\vskip2pt\hskip-1.5truecm{\sc{Bezdek, K.}} and {\sc{Connelly, R.}}

\bibitem[1991]{BezdekKConnelly91} Lower bounds for packing densities.
Acta Math. Hungar. {\bf{57}} (1991) no. 3-4, 291--311. MR1139324, 
DOI~10.1007/BF01903681

\vskip2pt\hskip-1.5truecm{\sc{Bezdek, K.}} and {\sc{L\'angi, Zs.}}

\bibitem[2020b]{BezdekKLangi20b} Bounds for totally separable translative packings in the plane.
Discrete Comput. Geom. {\bf{63}} (2020) no. 1, 49--72. MR4045741,
DOI~10.1007/s00454-018-0029-6

\bibitem[2022]{BezdekKLangi22} From the separable Tammes problem to extremal distributions of
great circles in the unit sphere. arXiv:2201.11234v1 [math.MG] 26 Jan 2022.

\vskip2pt\hskip-1.5truecm{\sc{Bilyk, D., Dai, F.}} and {\sc{Matzke, R. W.}}

\bibitem[2018]{BilykDaiMatzke} The Stolarsky principle and energy optimization on the sphere.
Constr. Approx. {\bf{48}} (2018) no. 1, 31--60. MR3825946,
DOI~10.1007/s00365-017-9412-4

\vskip2pt\hskip-1.5truecm{\sc{Bilyk, D.}} and {\sc{Matzke, R. W.}}

\bibitem[2019]{BilykMatzke} On the Fejes T\'oth problem about the sum of angles between lines.
Proc. Amer. Math. Soc. {\bf{147}} (2019) no. 1, 51--59. MR3876730,
DOI~10.1090/proc/14263

\vskip2pt\hskip-1.5truecm{\sc{Bleicher, M. N.}}

\bibitem[1975]{Bleicher75}  The thinnest three dimensional point lattice trapping a sphere.
Studia Sci. Math. Hungar. {\bf{10}} (1975) no. 1-2, 157--170. MR0476652

\vskip2pt\hskip-1.5truecm{\sc{Bleicher, M. N.}} and {\sc{Fejes T\'oth, L.}}

\bibitem[1964]{BleicherFTL64} Circle-packings and circle-coverings on a cylinder.
Michigan Math. J. {\bf{11}} (1964) 337--341. MR0169140, 
DOI~10.1307/mmj/1028999186

\vskip2pt\hskip-1.5truecm{\sc{Bleicher, M. N.}} and {\sc{Osborn, J. M.}}

\bibitem[1967]{BleicherOsborn} Minkowskian distributions of congruent discs.
Acta Math. Acad. Sci. Hungar. {\bf{18}} (1967) 5--17. MR0209978, 
DOI~10.1007/BF02020959

\vskip2pt\hskip-1.5truecm{\sc{Blind, G.}}

\bibitem[1972]{Blind72} Zug\"angliche Unterdeckungen der Ebene durch kongruente Kreise.
J. Reine Angew. Math. {\bf{257}} (1972) 29--46. MR0324551, 
DOI~10.1515/crll.1972.257.29

\bibitem[1976]{Blind76} $r$-zug\"angliche Unterdeckungen der Ebene durch kongruente Kreise. I.
J. Reine Angew. Math. {\bf{288}} (1976) 1--23. MR0442828, 
DOI~10.1515/crll.1976.288.1

\bibitem[1977]{Blind77} $r$-zug\"angliche Unterdeckungen der Ebene durch kongruente Kreise. II.
J. Reine Angew. Math. {\bf{289}} (1977) 1--29. MR0442829, 
DOI~0.1515/crll.1977.289.1

\bibitem[1981]{Blind81} $\varrho$-zug\"angliche Unterdeckungen der Sph\"are durch kongruente Kreise.
Resultate Math. {\bf{4}} (1981) no. 2, 141--154. MR0636482,
DOI~10.1007/BF03322974

\vskip2pt\hskip-1.5truecm{\sc{Blind, G.}} and {\sc{Blind, R.}}

\bibitem[1978]{BlindBlind78} Zug\"anglichkeit von Kugelpackungen im $R^n$.
Arch. Math. (Basel) {\bf{30}} (1978) no. 4, 438--439. MR0498196, 
DOI~10.1007/BF01226081

\bibitem[1979]{BlindBlind79} $r$-zug\"angliche Unterdeckungen der Ebene durch kongruente Kreise.
Studia Sci. Math. Hungar.  {\bf{14}}  (1979) no. 4, 441--452. MR0685911

\vskip2pt\hskip-1.5truecm{\sc{Bollob\'as, B.}}

\bibitem[1968a]{Bollobas68a} Remarks to a paper of L. Fejes T\'oth.
Studia Sci. Math. Hungar. {\bf{3}} (1968) 373--379. MR0239509 

\vskip2pt\hskip-1.5truecm{\sc{B\"or\"oczky, K.}}

\bibitem[1964]{Boroczky64} \"Uber stabile Kreis- und Kugelsysteme.
Ann. Univ. Sci. Budapest. E\"otv\"os Sect. Math. {\bf{7}} (1964) 79--82. MR0176391 

\bibitem[1967]{Boroczky67} \"Uber Dunkelwolken.
Proc. Colloquium on Convexity (Copenhagen, 1965) pp. 13--17
Kobenhavns Univ. Mat. Inst., Copenhagen. 1967. MR0215192 

\bibitem[1986]{Boroczky86} Closest packing and loosest covering of the space with balls.
Studia Sci. Math. Hungar.  {\bf{21}} (1986)  no. 1-2, 79--89. MR0898846 

\bibitem[2001]{Boroczky01} Edge close ball packings.
Discrete Comput. Geom. {\bf{26}} (2001) no. 1, 59--71. MR1832730,
DOI~10.1007/s00454-001-0021-3

\vskip2pt\hskip-1.5truecm{\sc{B\"or\"oczky, K.}}; {\sc{B\'ar\'any, I.}}; {\sc{Makai, E. Jr.}} and {\sc{Pach, J.}}

\bibitem[1986]{BoroczkyBaranyMakaiPach} Maximal volume enclosed by plates and proof of the chessboard conjecture.
Discrete Math. {\bf{60}} (1986) 101--120. MR0852101, 
DOI~10.1016/0012-365X(86)90006-3

\vskip2pt\hskip-1.5truecm{\sc{B\"or\"oczky, K.}} and {\sc{Soltan, V.}}

\bibitem[1996]{BoroczkySoltan} Translational and homothetic clouds for a convex body.
Studia Sci. Math. Hungar. {\bf{32}} (1996) no. 1-2, 93--102. MR1405128

\vskip2pt\hskip-1.5truecm{\sc{B\"or\"oczky, K.}} and {\sc{Szab\'o, L.}}

\bibitem[2002]{BoroczkySzabo02}  Minkowski arrangements of circles in the plane.
IV International Conference in "Stochastic Geometry, Convex Bodies, Empirical
Measures $\&$ Applications to Engineering Science", Vol. I (Tropea, 2001)
Rend. Circ. Mat. Palermo (2) Suppl.  No. 70, part I (2002) 87--92. MR1962558 

\bibitem[2004]{BoroczkySzabo04} Minkowski arrangements of spheres.
Monatsh. Math. {\bf{141}} (2004)  no. 1, 11--19. MR2109518, 
DOI~10.1007/s00605-002-0002-5

\vskip2pt\hskip-1.5truecm{\sc{B\"or\"oczky, K. Jr.}}

\bibitem[2004]{BoroczkyJr04} Finite packing and covering.
Cambridge Tracts in Mathematics, 154. Cambridge University Press, Cambridge, 2004. xviii+380 pp.
ISBN:0-521-80157-5. MR2078625

\vskip2pt\hskip-1.5truecm{\sc{B\"or\"oczky, K. Jr.}}; {\sc{F\'abi\'an, I.}} and {\sc{Wintsche, G.}}

\bibitem[2006]{BoroczkyJrFabianWintsche} Covering the crosspolytope by equal balls.
Period. Math. Hungar. {\bf{53}} (2006) no. 1-2, 103--113. MR2286463, 
DOI~10.1007/s10998-006-0024-1

\vskip2pt\hskip-1.5truecm{\sc{B\"or\"oczky, K., Jr.}} and {\sc{Tardos, G.}}

\bibitem[2002]{BoroczkyJrTardos} The longest segment in the complement of a packing.
Mathematika  {\bf{49}} (2002) no. 1-2, 45--49 MR2059040, 
DOI~10.1112/S002557930001603X

\vskip2pt\hskip-1.5truecm{\sc{B\"or\"oczky, K., Jr.}} and {\sc{Wintsche, G.}}

\bibitem[2000]{BoroczkyJrWintsche00} Sphere packings in the regular crosspolytope.
Ann. Univ. Sci. Budapest. E\"o\"otv\"os Sect. Math. {\bf{43}} (2000) 151--157 MR1847875

\vskip2pt\hskip-1.5truecm{\sc{Brandt, M.; Dickinson, W.; Ellsworth, A. V.; Kenkel, J.}} and {\sc{Smith, H.}}

\bibitem[2019]{BrandtDickinsonEllsworthKenkelSmith} Optimal packings of two to four equal circles on any flat torus.
Discrete Math. {\bf{342}} (2019) no. 12, 111597, 18 pp. MR3990010,
10.1016/j.disc.2019.111597

\vskip2pt\hskip-1.5truecm{\sc{Brauchart, J.~S.}} and {\sc{Dick, J.}}

\bibitem[2013]{BrauchartDick} A simple proof of Stolarsky's invariance principle.
Proc. Amer. Math. Soc. {\bf{141}} (2013) no. 6, 2085--2096. MR3034434,
DOI~10.1090/S0002-9939-2013-11490-5

\vskip2pt\hskip-1.5truecm{\sc{de Bruijn, N.G.}}

\bibitem[1954]{de Bruijn54} Problems 17 and 18 (in Dutch).
Nieuw Archief voor Wiskunde (3) {\bf{2}} (1954) 67.

\bibitem[1955]{de Bruijn55} Solution to problems 17 and 18 (in Dutch).
Wiskundige Opgaven met de oplossingen {\bf{20}} (1955) 19--20.

\vskip2pt\hskip-1.5truecm{\sc{Chan, Kwong-fai}} and {\sc{Lam, Tak Wah}}

\bibitem[1993]{ChanLam} An on-line algorithm for navigating in an unknown environment.
Selected papers from the 2nd Annual International Symposium on Algorithms, I (Taipei, 1991)
Internat. J. Comput. Geom. Appl. {\bf{3}} (1993) no. 3, 227--244. MR1241917, 
DOI~10.1142/S0218195993000154

\vskip2pt\hskip-1.5truecm{\sc{Chen, E. R.}}

\bibitem[2008]{ChenER} A dense packing of regular tetrahedra.
Discrete Comput. Geom. {\bf{40}} (2008) no. 2, 214--240. MR2438925,
DOI~10.1007/s00454-008-9101-y

\vskip2pt\hskip-1.5truecm{\sc{Chen, E. R.; Engel, M.}} and {\sc{Glotzer, S. C.}}

\bibitem[2010]{ChenEngelGlotzer} Dense crystalline dimer packings of regular tetrahedra.
Discrete Comput. Geom. {\bf{44}} (2010)  no. 2, 253--280. MR2671012,
DOI~10.1007/s00454-010-9273-0

\vskip2pt\hskip-1.5truecm{\sc{Connelly, R.}} and {\sc{Dickinson, W.}}

\bibitem[2014]{ConnellyDickinson} Periodic planar disc packings.
Philos. Trans. R. Soc. Lond. Ser. A Math. Phys. Eng. Sci. {\bf{372}}
(2014) no. 2008, 20120039, 17 pp. MR3158337,
DOI~10.1098/rsta.2012.0039

\vskip2pt\hskip-1.5truecm{\sc{Connelly, R., Funkhouser, M., Kuperberg, V.}} and {\sc{Solomonides, E.}}

\bibitem[2017]{ConnellyFunkhouserKuperbergSolomonides} Packings of equal disks in a square torus.
Discrete Comput. Geom. {\bf{58}} (2017) no. 3, 614--642. MR3690664,
DOI~10.1007/s00454-016-9843-x

\vskip2pt\hskip-1.5truecm{\sc{Connelly, R.; Shen, J. D.}} and {\sc{Smith, A. D.}}

\bibitem[2014]{ConnellyShenSmith} Ball packings with periodic constraints.
Discrete Comput. Geom. {\bf{52}} (2014) no. 4, 754--779. MR3279548,
DOI~10.1007/s00454-014-9636-z

\vskip2pt\hskip-1.5truecm{\sc{Conway, J. H.}} and {\sc{Torquato, S.}}

\bibitem[2006]{ConwayTorquato} Packing, tiling, and covering with tetrahedra.
Proc. Natl. Acad. Sci. USA {\bf{103}}  (2006)  no. 28, 10612--10617. MR2242647, 
DOI~10.1073/pnas.0601389103

\vskip2pt\hskip-1.5truecm{\sc{Cs\'oka, G.}}

\bibitem[1977]{Csoka} The number of congruent spheres that cover a given sphere of
three-dimensional space is not less than $30$. (Russian)
Studia Sci. Math. Hungar. {\bf{12}} (1977) no. 3-4, 323--334. MR0607086 

\vskip2pt\hskip-1.5truecm{\sc{Danzer, L.}}

\bibitem[1960]{Danzer60} Drei Beispiele zu Lagerungsproblemen.
Arch. Math. (Basel)  {\bf{11}} (1960) no. 1, 159--165. MR1552702,
DOI~10.1007/BF01236925

\vskip2pt\hskip-1.5truecm {\sc{Davenport, H.}} and {\sc{Haj\'os, G.}}

\bibitem[1951]{DavenportHajos} Problem 35.
Matematikai Lapok {\bf{2}} (1951) 63.

\vskip2pt\hskip-1.5truecm{\sc{Dawson, R.}}

\bibitem[1984]{Dawson} On removing a ball without disturbing the others.
Math. Mag. {\bf{57}} (1984) no. 1, 27--30. MR0729036, 
DOI~10.1080/0025570X.1984.11977071

\vskip2pt\hskip-1.5truecm{\sc{Delone, B. N.}} and {\sc{Ry\v skov, S. S.}}

\bibitem[1963]{DeloneRyskov} Solution of the problem on the least dense lattice
covering of a 4-dimensional space by equal spheres. (Russian)
Dokl. Akad. Nauk SSSR {\bf{152}} (1963) 523--524.
English translation in Soviet Math. Dokl. {\bf{4}} (1963) 1333--1334. MR0175850

\vskip2pt\hskip-1.5truecm{\sc{Dickinson, W.; Guillot, D.; Keaton, A.}} and {\sc{Xhumari, S.}}

\bibitem[2011a]{DickinsonGuillotKeatonXhumari11a} Optimal packings of up to
six equal circles on a triangular flat torus.
J. Geom. {\bf{102}} (2011) no. 1-2, 27--51. MR2904615,
DOI~10.1007/s00022-011-0099-6

\bibitem[2011b]{DickinsonGuillotKeatonXhumari11b} Optimal packings of up to
five equal circles on a square flat torus.
Beitr\"age Algebra Geom. {\bf{52}} (2011) no. 2, 315--333. MR2842632

\vskip2pt\hskip-1.5truecm{\sc{Dominy\'ak, I.}}

\bibitem[1964]{Dominyak} On the density of stable systems of circles. (Hungarian)
Magyar Tud. Akad. Mat. Fiz. Oszt. K\"ozl. {\bf{14}} (1964) 401--413. MR0229141 

\vskip2pt\hskip-1.5truecm{\sc{Dougherty, R.}} and {\sc{Faber, V.}}

\bibitem[2004]{DoughertyFaber} The degree-diameter problem for several varieties of Cayley graphs. I. The abelian case.
SIAM J. Discrete Math. {\bf{17}} (2004) 478--519. MR2050686, 
DOI~10.1137/S0895480100372899

\vskip2pt\hskip-1.5truecm{\sc{Dumir, V. C.}} and {\sc{Khassa, D. S.}}

\bibitem[1973a]{DumirKhassa73a} A conjecture of Fejes T\'oth on saturated systems of circles.
Proc. Cambridge Philos. Soc. {\bf{74}} (1973) 453--460. MR0333998, 
DOI~10.1017/S0305004100077203

\bibitem[1973b]{DumirKhassa73b} Saturated systems of symmetric convex domains; results of Eggleston, Bambah and Woods.
Proc. Cambridge Philos. Soc. {\bf{74}} (1973) 107--116. MR0315599, 
DOI~10.1017/S030500410004785X

\vskip2pt\hskip-1.5truecm{\sc{Eggleston, H. G.}}

\bibitem[1965]{Eggleston65} A minimal density plane covering problem.
Mathematika {\bf{12}} (1965) 226--234. MR0192413, 
DOI~10.1112/S0025579300005325

\vskip2pt\hskip-1.5truecm{\sc{Fejes T\'oth, G.}}

\bibitem[1976b]{FTG76b} Research problem no. 18.
Periodica Math. Hungar. {\bf{7}} (1976) no. 1,  89--90. MR1553596 DML,
DOI~10.1007/BF02019998

\bibitem[1978]{FTG78} Evading convex discs.
Studia Sci. Math. Hungar. {\bf{13}} (1978) no. 3-4, 453--461. MR0620156 

\bibitem[1987]{FTG87} Totally separable packing and covering with circles.
Studia Sci. Math. Hungar. {\bf{22}} (1987) no. 1-4, 65--73. MR0913893 

\bibitem[2013]{FTG13} Shortest path avoiding circles and balls.
Studia Sci. Math. Hungar. {\bf{50}} (2013) 454--464. MR3187827

\vskip2pt\hskip-1.5truecm{\sc{Fejes T\'oth, G.}} and {\sc{Fejes T\'oth, L.}}

\bibitem[1973b]{FTGFTL73b} On totally separable domains.
Acta Math. Acad. Sci. Hungar. {\bf{24}} (1973) 229--232. MR0322690, 
DOI~10.1007/BF01894631

\vskip2pt\hskip-1.5truecm{\sc{Fejes T\'oth, G.; Kuperberg, G.}} and {\sc{Kuperberg, W.}}

\bibitem[1998]{FTGKuperbergKuperberg} Highly saturated packings and reduced coverings.
Monatsh. Math. {\bf{125}}  (1998)  no. 2, 127--145. MR1604938, 
DOI~10.1007/BF01332823

\vskip2pt\hskip-1.5truecm{\sc{Fejes T\'oth, L.}}

\bibitem[1956b]{FTL56b} On the sum of distances determined by a pointset.
Acta Math. Acad. Sci. Hungar {\bf{7}} (1956) 397--401. MR0107212, 
DOI~10.1007/BF02020534

\bibitem[1959b]{FTL59b} An extremal distribution of great circles on a sphere.
Publ. Math. Debrecen {\bf{6}} (1959) 79--82. MR0105065 

\bibitem[1959d]{FTL59d} \"Uber eine Punktverteilung auf der Kugel.
Acta Math. Acad. Sci. Hungar. {\bf{10}} (1959) 13--19. MR0105654, 
DOI~10.1007/BF02063286

\bibitem[1959e]{FTL59e} Verdeckung einer Kugel durch Kugeln.
Publ. Math. Debrecen {\bf{6}} (1959) 234--240. MR0113179 

\bibitem[1960a]{FTL60a} On the stability of a circle packing.
Ann. Univ. Sci. Budapest. E\"otv\"os Sect. Math. {\bf{3-4}} (1960/1961) 63--66. MR0140005 

\bibitem[1962]{FTL62} Dichteste Kreispackungen auf einem Zylinder.
Elem. Math. {\bf{17}} (1962) 30--33. MR0133738 

\bibitem[1965b]{FTL65b} Minkowskian distribution of discs.
Proc. Amer. Math. Soc. {\bf{16}} (1965) 999--1004. MR0180921, 
DOI~0.1090/S0002-9939-1965-0180921-7

\bibitem[1966b]{FTL66b} On the permeability of a circle-layer.
Studia Sci. Math. Hungar. {\bf{1}} (1966) 5--10. MR0206824 

\bibitem[1967a]{FTL67a} Packings and coverings in the plane.
1967 Proc. Colloquium on Convexity (Copenhagen, 1965) pp. 78-87
Kobenhavns Univ. Mat. Inst., Copenhagen. MR0217702 

\bibitem[1967b]{FTL67b} Minkowskian circle-aggregates.
Math. Ann. {\bf{171}} (1967) 97--103. MR0221386, 
DOI~10.1007/BF01351644

\bibitem[1967c]{FTL67c} On the arrangement of houses in a housing estate.
Studia Sci. Math. Hungar. {\bf{2}} (1967) 37--42. MR0215188 

\bibitem[1968a]{FTL68a} On the permeability of a layer of parallelograms.
Studia Sci. Math. Hungar. {\bf{3}} (1968) 195--200. MR0232285 

\bibitem[1973d]{FTL73d} On the density of a connected lattice of convex bodies.
Acta Math. Acad. Sci. Hungar. {\bf{24}} (1973) 373--376. MR0331228,
DOI~10.1007/BF01958049

\bibitem[1973e]{FTL73e} Distribution of points on convex polyhedra.
Period. Math. Hungar. {\bf{4}} (1973) 29--38. MR0331226,
DOI~10.1007/BF02018034

\bibitem[1975b]{FTL75b} Research problem no. 13.
Period. Math. Hungar. {\bf{6}} (1975) no. 2, 197--199. MR1553592 DML, 
DOI~10.1007/BF02018822

\bibitem[1975c]{FTL75c} Research problem no. 14.
Period. Math. Hungar. {\bf{6}} (1975) no. 3, 277--278. MR1553594 DML, 
DOI~10.1007/BF02018282

\bibitem[1976]{FTL76} Close packing and loose covering with balls.
Publ. Math. Debrecen {\bf{23}} (1976) no. 3-4, 323--326. MR0428199 

\bibitem[1978a]{FTL78a} Remarks on the closest packing of convex discs.
Comment. Math. Helv. {\bf{53}} (1978) no. 4, 536--541. MR0514025, 
DOI~10.1007/BF02566097

\bibitem[1978b]{FTL78b} Research problem no. 24.
Periodica Math. Hungar. {\bf{9}} (1978) no. 1-2, 173--174. MR1553604 DML,
DOI~10.1007/BF02018932

\bibitem[1987]{FTL87} On spherical tilings generated by great circles.
Geom. Dedicata {\bf{23}} (1987) no. 1, 67--71. MR0886775, 
DOI~10.1007/BF00147392

\bibitem[1984b]{FTL84b} Compact packing of circles.
Studia Sci. Math. Hungar. {\bf{19}} (1984) no. 1, 103--107. MR0787791 

\bibitem[1993]{FTL93} Flight in a packing of disks.
Discrete Comput. Geom.  {\bf{9}}  (1993)  no. 1, 1--9. MR1184690, 
DOI~10.1007/BF02189303

\bibitem[1999]{FTL99} Minkowski circle packings on the sphere.
Discrete Comput. Geom. {\bf{22}}  (1999)  no. 2, 161--166. MR1698538, 
DOI~10.1007/PL00009451

\vskip2pt\hskip-1.5truecm{\sc{Fejes T\'oth, L.}} and {\sc{Heppes, A.}}

\bibitem[1963]{FTLHeppes63} \"Uber stabile K\"orpersysteme.
Compositio Math. {\bf{15}} (1963) 119--126. MR0161227 

\bibitem[1980]{FTLHeppes80} Multi-saturated packings of circles.
Studia Sci. Math. Hungar. {\bf{15}} (1980) no. 1-3, 303--307. MR0681452

\vskip2pt\hskip-1.5truecm{\sc{Fejes T\'oth, L.}} and {\sc{Makai, E., Jr.}}

\bibitem[1974]{FTLMakai} On the thinnest non-separable lattice of convex plates.
Studia Sci. Math. Hungar.  {\bf{9}}  (1974) 191--193. MR0370369

\vskip2pt\hskip-1.5truecm{\sc{Fernique, T..}}

\bibitem[2019a]{Fernique19a} A densest ternary circle packing in the plane.
arXiv:1912.02297v1 [cs.CG] 4 Dec 2019


\bibitem[2019b]{Fernique19b} Compact packings of space with three sizes of spheres.
arXiv:1912.02293v1 [math.MG] 4 Dec 2019

\bibitem[2021]{Fernique21} Compact packings of space with two sizes of spheres.
Discrete Comput. Geometry  {\bf 65} (2021) no.\ 4, 1287--1295. MR4249904,
DOI~10.1007/s00454-019-00140-8

\vskip2pt\hskip-1.5truecm{\sc{Fernique, T.; Hashemi, A.}} and {\sc{Sizova, O.}}


\bibitem[2020]{FerniqueHashemiSizova} Compact packings of the plane with three sizes of discs.
Discrete Comput. Geometry {\bf 66} (2021) no.\ 2, 613--635. MR4292755,
DOI~10.1007/s00454-019-00166-y

\vskip2pt\hskip-1.5truecm{\sc{Fernique, T}} and {\sc{Pchelina, D.}}

\bibitem[2021]{FerniquePchelina}
Compact packings are not always the densest.
arXiv:2104.12458v2 [math.MG] 27 Apr 2021

\vskip2pt\hskip-1.5truecm{\sc{Fiduccia, C. M. Forcade, R. W.}} and {\sc{Zito, J. S.}}

\bibitem[1998]{FiducciaForcadeZito} Geometry and diameter bounds of directed Cayley graphs of abelian groups.
SIAM J. Discrete Math. {\bf{11}} (1998) 157--167. MR1612881, 
DOI~10.1137/S0895480195286456

\vskip2pt\hskip-1.5truecm{\sc{Florian, A.}}

\bibitem[1967]{Florian67} Zur Geometrie der Kreislagerungen.
Acta Math. Acad. Sci. Hungar. {\bf{18}} (1967) 341--358. MR0218980, 
DOI~10.1007/BF02280295

\bibitem[1968]{Florian68} Bemerkung zu einer Arbeit von L. Fejes T\'oth.
Studia Sci. Math. Hungar. {\bf{3}} (1968) 359--363. MR0236821 

\bibitem[1978b]{Florian78b} On the permeability of layers of discs.
Studia Sci. Math. Hungar. {\bf{13}} (1978) no. 1-2, 125--132. MR0630383 

\bibitem[1979b]{Florian79b} \"Uber die Durchl\"assigkeit gewisser Scheiben\-schich\-ten.
\"Oster\-reich. Akad. Wiss. Math.-Natur. Kl. Sitzungsber. II
{\bf{188}}  (1979) no. 8-10, 417--427. MR0599882 

\bibitem[1980a]{Florian80a} \"Uber die Durchl\"assigkeit einer Schicht konvexer Scheiben.
Studia Sci. Math. Hungar.  {\bf{15}}  (1980) no. 1-3, 201--213.
MR0681440 

\bibitem[1980b]{Florian80b} Durchl\"assigkeit einer Scheibenschicht.
Twelfth Styrian Mathematical Symposium (Graz, 1980), Ber. No. 144, 17 pp., Bericht,
140--150, Forschungszentrum Graz, Graz, 1980. MR0625870

\bibitem[1985]{Florian85} On compact packing of circles.
Studia Sci. Math. Hungar. {\bf{20}} (1985) no. 1-4, 473--480. MR0886051 

\vskip2pt\hskip-1.5truecm{\sc{Florian, A.}} and {\sc{Groemer, H.}}

\bibitem[1985]{FlorianGroemer} Two remarks on the permeability of layers of convex bodies.
Studia Sci. Math. Hungar. {\bf{20}}  (1985)  no. 1-4, 259--265. MR0886028 

\vskip2pt\hskip-1.5truecm{\sc{Fodor, F.; V\'{\i}gh, V.}} and {\sc{Zarn\'ocz, T.}}

\bibitem[2016a]{FodorVighZarnocz16a} On the angle sum of lines.
Arch. Math. (Basel) {\bf{106}} (2016) no. 1, 91--100. MR3451371,
DOI~10.1007/s00013-015-0847-1

\vskip2pt\hskip-1.5truecm{\sc{F\"oldv\'ri, V.}}

\bibitem[2020]{Foldvari} Bounds on convex bodies in pairwise intersecting Minkowski arrangement of order $\mu$.
J. Geom. (2020) 111:27
DOI~10.1007/s00022-020-00538-3

\vskip2pt\hskip-1.5truecm{\sc{Forcade, R.}} and {\sc{Lamoreaux, J.}}

\bibitem[2000]{ForcadeLamoreaux} Lattice-simplex coverings and the 84-shape.
SIAM J. Discrete Math. {\bf{13}} (2000) 194--201. MR1760337, 
DOI~10.1137/S0895480198349622

\vskip2pt\hskip-1.5truecm{\sc{Frostman, O.}}

\bibitem[1953]{Frostman} A theorem of F\'ary with elementary applicatons. (Swedish)
Nordisk Mat. Tidskr. {\bf{1}} (1953) 25--32, 64. MR0055696 

\vskip2pt\hskip-1.5truecm{\sc{F\"uredi, Z.}}

\bibitem[1991]{Furedi91} The densest packing of equal circles into a parallel strip.
Discrete Comput. Geom. {\bf{6}} (1991) no. 2, 95--106. MR1083626,
DOI~10.1007/BF02574677

\vskip2pt\hskip-1.5truecm{\sc{F\"uredi, Z.}} and {\sc{Loeb, P. A.}}

\bibitem[1994]{FurediLoeb} On the best constant for the Besicovitch covering theorem.
Proc. Amer. Math. Soc. {\bf{121}} (1994) no. 4, 1063--1073. MR1249875,
DOI~10.1090/S0002-9939-1994-1249875-4

\vskip2pt\hskip-1.5truecm{\sc{Goldberg, M.}}

\bibitem[1971]{Goldberg71} On the densest packing of equal spheres in a cube.
Math. Mag. {\bf{44}} (1971) 199--208. MR0298562,
DOI~10.1080/0025570X.1971.11976147

\vskip2pt\hskip-1.5truecm{\sc{Golser, G.}}

\bibitem[1977]{Golser} Dichteste Kugelpackungen im Oktaeder.
Studia Sci. Math. Hungar. {\bf{12}} (1977) no. 3-4, 337--343. MR0607088

\vskip2pt\hskip-1.5truecm{\sc{Gonz\'alez Merino, B.}} and {\sc{Schymura, M.}}

\bibitem[2017]{GonzalezMerinoSchymura} On densities of lattice arrangements
intersecting every $i$-dimensional affine subspace.
Discrete Comput. Geom. {\bf{58}} (2017) no. 3, 663--685. MR3690667,
DOI~10.1007/s00454-017-9911-x

\vskip2pt\hskip-1.5truecm{\sc{Graf, C.}} and {\sc{Paukowitsch, P.}}

\bibitem[1997]{GrafPaukowitsch} M\"oglichst dichte Packungen aus kongruenten
Drehzylindern mit paarweise windschiefen Achsen.
Elem. Math. {\bf{52}} (1997) 71--83. MR1450459, 
DOI~10.1007/s000170050013

\vskip2pt\hskip-1.5truecm{\sc{Gravel, S.}}; {\sc{Elser, V.}} and {\sc{Kallus, Y.}}

\bibitem[2011]{GravelElserKallus} Upper bound on the packing density of regular tetrahedra and octahedra.
Discrete Comput. Geom. {\bf{46}} (2011) no. 4, 799--818. MR2846180,
DOI~10.1007/s00454-010-9304-x

\vskip2pt\hskip-1.5truecm{\sc{Gritzmann, P.}} and {\sc{Wills, J. M.}}

\bibitem[1993]{GritzmannWills} Finite packing and covering.
Handbook of convex geometry, Vol. A, B, 861--897, North-Holland, Amsterdam, 1993. MR1242998,
DOI~10.1016/B978-0-444-89597-4.50008-1

\vskip2pt\hskip-1.5truecm{\sc{Groemer, H.}}

\bibitem[1960a]{Groemer60a} \"Uber die Einlagerung von Kreisen in einen konvexen Bereich.
Math. Z. {\bf{73}} (1960) 285--294. MR0110981, 
DOI~10.1007/BF01159721

\bibitem[1966b]{Groemer66b} Zusammenh\"angende Lagerungen konvexer K\"orper.
Math. Z. {\bf{94}} (1966) 66--78. MR0199792, 
DOI~10.1007/BF01111261

\vskip2pt\hskip-1.5truecm{\sc{Gr\"unbaum, B.}}

\bibitem[1960]{Grunbaum60} On a problem of L. Fejes-T\'oth.
Amer. Math. Monthly {\bf{67}} (1960) 882--884. MR0139974, 
DOI~10.1080/00029890.1960.11992016

\vskip2pt\hskip-1.5truecm{\sc{Haji-Akbari, A.; Engel, M., Keys, A.}}; {\sc{Zheng, X.}}; {\sc{Petschek, R. G.}}; {\sc{Palffy-Muhoray, P.}} and
{\sc{Glotzer S. C.}}

\bibitem[2009]{Haji-Akbari+} Disordered, quasicrystalline and crystalline phases of densely packed tetrahedra.
Nature {\bf{462}} (2009) 773--777.
DOI~10.1038/nature08641

\vskip2pt\hskip-1.5truecm{\sc{Haj\'os, G.}}

\bibitem[1964]{Hajos64} \"Uber Kreiswolken.
Ann. Univ. Sci. Budapest. E\"otv\"os Sect. Math. {\bf{7}} (1964) 55--57. MR0175041 

\vskip2pt\hskip-1.5truecm{\sc{Hausel, T.}}

\bibitem[1992]{Hausel} Transillumination of lattice packing of balls.
Studia Sci. Math. Hungar. {\bf{27}} (1992) no. 1-2, 241--242. MR1207577 

\vskip2pt\hskip-1.5truecm{\sc{Henk, M.}}

\bibitem[2005]{Henk} Free planes in lattice sphere packings.
Adv. Geom. {\bf{5}} (2005) no. 1, 137--144. MR2110466, 
DOI~10.1515/advg.2005.5.1.137

\vskip2pt\hskip-1.5truecm{\sc{Henk, M}}; {\sc{Ziegler, G. M.}} and {\sc{Zong, C.}}

\bibitem[2002]{HenkZieglerZong} On free planes in lattice ball packings.
Bull. London Math. Soc.  {\bf{34}}  (2002)  no. 3, 284--290. MR1887700, 
DOI~10.1112/S0024609301008888

\vskip2pt\hskip-1.5truecm{\sc{Henk, M.}} and {\sc{Zong, C.}}

\bibitem[2000]{HenkZong00} Segments in ball packings.
Mathematika {\bf{47}} (2000) no. 1-2, 31--38. MR1924485,
DOI~10.1112/S0025579300015692

\vskip2pt\hskip-1.5truecm{\sc{Heppes, A.}}

\bibitem[1958]{Heppes58} An extremal property of the spherical net of the cuboctahedron. (Hungarian)
Magyar Tud. Akad. Mat. Kutat\'o Int. K\"ozl. {\bf{8}} (1958) no. 1/2, 97--99. MR0107216 

\bibitem[1960b]{Heppes60b} Ein Satz \"uber gitterf\"ormige Kugelpackungen.
Ann. Univ. Sci. Budapest. E\"otv\"os Sect. Math. {\bf{3--4}} (1960/1961) 89--90. MR0133737 

\bibitem[1961a]{Heppes61a} \"Uber Kreis- und Kugelwolken.
Acta Math. Acad. Sci. Hungar. {\bf{12}} (1961) 209--214. MR0126213, 
DOI~10.1007/BF02066683

\bibitem[1967a]{Heppes67a} On the densest packing of circles not blocking each other.
Studia Sci. Math. Hungar. {\bf{2}} (1967) 257--263. MR0215194 

\bibitem[1967b]{Heppes67b} On the number of spheres which can hide a given sphere.
Canad. J. Math. {\bf{19}} (1967) 413--418. MR0209980, 
DOI~10.4153/CJM-1967-033-4

\bibitem[1999]{Heppes99} Densest circle packing on the flat torus.
Periodica Math. Hunar. {\bf{39}} (1999) no. 1-3, 129--134. MR1783819,
DOI~10.1023/A:1004847008346

\bibitem[2000]{Heppes00} On the densest packing of discs of radius $1$ and $\sqrt2-1$.
Studia Sci. Math. Hungar. {\bf{36}} (2000) no. 3-4, 433--454. MR1798749

\bibitem[2001b]{Heppes01b} On the density of 2-saturated lattice packings of discs.
Monatsh. Math. {\bf{134}} (2001) no. 1, 51--66. MR1872046, 
DOI~10.1007/s006050170011

\bibitem[2002]{Heppes02} Expandability radius versus density of a lattice packing.
Period. Math. Hungar. {\bf{45}} (2002) no. 1-2, 65--71. MR1955193,
DOI~10.1023/A:1022345929586

\bibitem[2003b]{Heppes03b} Some densest two-size disc packings in the plane.
Disicrete Comput.~Geom. {\bf{30}} (2003) 241--262. MR2007963,
DOI~10.1007/s00454-003-0007-6

\vskip2pt\hskip-1.5truecm{\sc{Heppes, A.}} and {\sc{Kert\'esz, G.}}

\bibitem[1997]{HeppesKertesz} Packing circles of two different sizes on the sphere.
Intuitive geometry (Budapest, 1995) Bolyai Soc. Math. Stud., 6,
J\'anos Bolyai Math. Soc., Budapest, 1997, 357--365. MR1470773

\vskip2pt\hskip-1.5truecm{\sc{Hilbert, D.}}

\bibitem[1900]{Hilbert} Mathematical problems.
Reprinted from Bull. Amer. Math. Soc. {\bf{8}} (1902) 437--479.
Bull. Amer. Math. Soc. (N.S.) {\bf{37}} (2000) no. 4, 407--436.
Earlier publications (in the original German) appeared in G\"ottinger Nachrichten, 1900, 253--297,
and Archiv der Mathematik und Physik, (3) {\bf{1}} (1901) 44--63, 213--237. MR1779412,
DOI~10.1090/S0002-9904-1902-00923-3

\vskip2pt\hskip-1.5truecm{\sc{Hortob\'agyi, I.}}

\bibitem[1971]{Hortobagyi71} Durchleuchtung gitterf\"ormiger Kugelpackungen mit Lichtb\"undeln.
Studia Sci. Math. Hungar. {\bf{6}} (1971) 147--150. MR0348634 

\bibitem[1976a]{Hortobagyi76a} \"Uber die Durchl\"assigkeit einer
aus Scheiben konstanter Breite bestehenden Schicht.
Studia Sci. Math. Hungar. {\bf{11}} (1976) no. 3-4, 383--387. MR0554584 

\vskip2pt\hskip-1.5truecm{\sc{Horv\'ath, J.}}

\bibitem[1974]{Horvath74} Die Dichte einer Kugelpackung in einer 4-dimensionalen Schicht.
Period. Math. Hungar. 5 (1974) 195--199. MR0376543
DOI 10.1007/BF02023199

\bibitem[1980]{Horvath80} Narrow latticed packing of unit balls in the space $E^n$, (Russian)
in Geometry of Positive Quadratic Forms, Trudy Mat. Inst. Steklov. {\bf{152}} (1980) 216--231, 238.
English translation in Proc. Steklov Inst. Math. {\bf{152}} no. 3, (1982) 237--254.  MR0603826 

\bibitem[1986]{Horvath86} Several problems of $n$-dimensional discrete geometry. (Russian)
Doctoral thesis, Steklov Institute of Mathematics, 1986.

\vskip2pt\hskip-1.5truecm{\sc{Horv\'ath, J.}} and {\sc{Moln\'ar, J.}}

\bibitem[1967]{HorvathMolnar67} On the density of non-overlapping unit spheres lying in a strip.
Ann. Univ. Sci. Budapest. E\"otv\"os Sect. Math. {\bf{10}} (1967) 193--201. MR0243431 

\vskip2pt\hskip-1.5truecm{\sc{Horv\'ath, J.}} and {\sc{Ry{\u{s}}kov, S. S.}}

\bibitem[1975a]{HorvathRyskov75a} On the radii of the cylinders which fit
into a lattice packing of the $n$-dimensional space with unit spheres. (Hungarian)
Mat. Lapok {\bf{26}} (1975) no. 1-2, 91--96. MR0514022 

\bibitem[1975b]{HorvathRyskov75b} Estimation of the radius of a cylinder that
can be imbedded in every lattice packing of $n$-dimensional unit balls. (Russian)
Mat. Zametki {\bf{17}} (1975) 123--128.
English translation: Math. Notes 17 (1975), nos. 1--2, 72--75. MR0370373,
DOI~10.1007/BF01093847

\vskip2pt\hskip-1.5truecm{\sc{Ismailescu, D.}}

\bibitem[2003]{Ismailescu03} Slicing the pie.
U.S.-Hungarian Workshops on Discrete Geometry and Convexity (Budapest, 1999/ Auburn, AL, 2000).
Discrete Comput. Geom. {\bf{30}} (2003) no. 2, 263--276. MR2007964, 
DOI~10.1007/s00454-003-0008-5

\vskip2pt\hskip-1.5truecm{\sc{Ismailescu, D.}} and {\sc{Laskawiec, P.}}

\bibitem[2019]{IsmailescuLaskawiec} Dense packings with nonparallel cylinders.
Elem. Math. {\bf{74}} (2019) no. 3, 89--103. MR3981287,
DOI~10.4171/EM/387

\vskip2pt\hskip-1.5truecm{\sc{Jo\'os, A.}}

\bibitem[2008b]{Joos08b} Covering the unit cube by equal balls.
Beitr\"age Algebra Geom. {\bf{49}} (2008) no. 2, 599--605. MR2468076

\bibitem[2009a]{Joos09a} On the packing of fourteen congruent spheres in a cube.
Geom. Dedicata {\bf{140}} (2009) 49--80. MR2504734

\bibitem[2009b]{Joos09b} Erratum to: A. Jo\'os: Covering the unit cube by equal balls [MR2468076].
Beitr\"age Algebra Geom. {\bf{50}} (2009) no. 2, 603--605. MR2572024

\bibitem[2014a]{Joos14a} Covering the $k$-skeleton of the 3-dimensional unit cube by five balls.
Beitr\"age Algebra Geom. {\bf{55}} (2014) no. 2, 393--414. MR3263252,
DOI~10.1007/s13366-013-0144-8

\bibitem[2014b]{Joos14b} Covering the $k$-skeleton of the 3-dimensional unit cube by six balls.
Discrete Math. {\bf{336}} (2014) 85--95. MR3254974,
DOI~10.1016/j.disc.2014.07.017

\bibitem[2018]{Joos18} Covering the 5-dimensional unit cube by eight congruent balls.
Period. Math. Hungar. {\bf{77}} (2018) no. 1, 77--82. MR3842971,
DOI~10.1007/s10998-018-0241-4

\bibitem[2019]{Joos19} On covering the square flat torus by congruent discs.
Australas. J. Combin. {\bf{75}} (2019) 113--126. MR3997953

\vskip2pt\hskip-1.5truecm{\sc{Jo\'os, A.}} and {\sc{Nagy, B.}}

\bibitem[2008a]{JoosNagy} Optimal packings of 2, 3, and 4 equal balls into a cubical flat 3-torus.
Boll. Unione Mat. Ital. {\bf{13}} (2020) no. 3, 335--340. MR4132911,
DOI~10.1007/s40574-020-00224-x

\vskip2pt\hskip-1.5truecm{\sc{Jucovi\v{c}, E.}}

\bibitem[1966]{Jucovic66} \"Uber die minimale Dicke einer $k$-fachen Kreiswolke.
Ann. Univ. Sci. Budapest. E\"otv\"os Sect. Math. {\bf{9}} (1966) 143--146. MR0206825 

\bibitem[1970]{Jucovic70} Raumanspr\"uchliche Kreispackungen in der euklidischen Ebene.
Mat. \v{C}asopis Sloven. Akad. Vied {\bf{20}} (1970) 3--10. MR0310769 

\vskip2pt\hskip-1.5truecm {\sc{Kadlicsk\'o, M.}} and {\sc{L\'angi, Zs.}}

\bibitem[2022]{KadlicskoLangi} On generalized Minkowski arrangements.
Ars Mathematica Contemporanea (2022)
DOI~10.26493/1855-3974.2550.d96

\vskip2pt\hskip-1.5truecm {\sc{Kahle, M.}}

\bibitem[2012]{Kahle} Sparse locally-jammed disk packings.
Ann. Comb. {\bf{16}} (2012) no. 4, 773--780. MR3000444,
DOI~10.1007/s00026-012-0159-0

\hskip-1.5truecm{\sc Kallus, Y.}; {\sc Elser, V.} and {\sc Gravel, S.}

\bibitem[2010]{KallusElserGravel} Dense periodic packings of tetrahedra with small repeating units.
Discrete Comput. Geom. {\bf 44} (2010) no. 2, 245--252. MR267101,
DOI~10.1007/s00454-010-9254-3

\vskip2pt\hskip-1.5truecm{\sc Kannan, R.} and {\sc Lov\'{a}sz, L.}

\bibitem[1988]{KannanLovasz} Covering minima and lattice-point-free convex bodies.
Annals of Mathematics (2) {\bf{128}} (1988) 577-602. MR0970611,
DOI~10.2307/1971436

\vskip2pt\hskip-1.5truecm{\sc{Kennedy, T.}}

\bibitem[2004]{Kennedy04} A densest compact planar packing with two sizes of discs.
arXiv:math/0412418v1 [math.MG] 21 Dec 2004

\bibitem[2006]{Kennedy06} Compact packings of the plane with two sizes of discs.
Discrete Comput. Geom.  {\bf{35}}  (2006)  no. 2, 255--267. MR2195054, 
DOI~10.1007/s00454-005-1172-4

\vskip2pt\hskip-1.5truecm{\sc{Kert\'esz, G.}}

\bibitem[1982]{Kertesz82} Thesis, E\"otv\"os Lor\'and University of Budapest, 1982.

\bibitem[1988]{Kertesz88} On totally separable packings of equal balls.
Acta Math. Hungar.  {\bf{51}}  (1988)  no. 3-4, 363--364. MR0956988, 
DOI~10.1007/BF01903343

\vskip2pt\hskip-1.5truecm{\sc{Khassa, D. S.}}

\bibitem[1975a]{Khassa75a} A lower bound for the lower densities of saturated systems of convex domains.
Indian J. Pure Appl. Math. {\bf{6}} (1975) no. 12, 1422--1435. MR0487803

\bibitem[1975b]{Khassa75b} A lower bound for the lower densities of saturated systems of spheres.
Indian J. Pure Appl. Math. {\bf{6}} (1975) no. 12, 1436--1440. MR0487804

\vskip2pt\hskip-1.5truecm{\sc{Kuperberg, K.}}

\bibitem[1990]{KuperbergK} A nonparallel cylinder packing with positive density.
Mathematika {\bf{37}} (1990) 324--331. MR1099780, 
DOI~10.1112/S0025579300013036

\vskip2pt\hskip-1.5truecm{\sc{Kuperberg, W.}}

\bibitem[2007]{Kuperberg07} Optimal arrangements in packing congruent balls in a spherical container.
Discrete Comput. Geom. {\bf{37}} (2007) no. 2, 205--212. MR2295053, 
DOI~10.1007/s00454-006-1303-6

\vskip2pt\hskip-1.5truecm{\sc{Lagarias, J. C.}} and {\sc{Zong, C.}}

\bibitem[2012]{LagariasZong} Mysteries in packing regular tetrahedra.
Notices of the AMS {\bf{58}} (2012) 1540--1549. MR3027108

\vskip2pt\hskip-1.5truecm{\sc{Larcher, H.}}

\bibitem[1962]{Larcher} Solution of a geometric problem by Fejes T\'oth.
Michigan Math. J. {\bf{9}} (1962) 45--51. MR0137032, 
DOI~10.1307/mmj/1028998619

\vskip2pt\hskip-1.5truecm{\sc{Linhart, J.}}

\bibitem[1978]{Linhart78} Closest packings and closest coverings by translates of a convex disc.
Studia Sci. Math. Hungar. {\bf{13}} (1978) no. 1-2, 157--162. MR0630388 

\vskip2pt\hskip-1.5truecm{\sc{Maehara, H.}}

\bibitem[1995]{Maehara95} An extremal problem for arrangement of great circles.
Math. Japon. {\bf{41}} (1995) 125--129. MR1317755

\vskip2pt\hskip-1.5truecm{\sc{Makai, E., Jr.}}

\bibitem[1978]{Makai78} On the thinnest nonseparable lattice of convex bodies.
Studia Sci. Math. Hungar. {\bf{13}} (1978) no. 1-2, 19--27. MR0630376 

\vskip2pt\hskip-1.5truecm{\sc{Makai, E., Jr.}}. and {\sc{Martini, H.}}

\bibitem[2016]{MakaiMartini} Density estimates for $k$-impassable lattices of balls and general convex bodies in $R^n$.
arXiv:1612.01307v1 [math.MG] 5 Dec 2016

\vskip2pt\hskip-1.5truecm{\sc{Melissen, H}}

\bibitem[1997]{Melissen97} Packing and covering with circles.
Thesis, University Utrecht 1997. ISBN 90-393-1500-0

\vskip2pt\hskip-1.5truecm{\sc{Messerschmidt, M.}}

\bibitem[2020]{Messerschmidt20} On compact packings of the plane with circles of three radii.
Comput. Geom. {\bf{86}} (2020) 101564, 17 pp. MR4033125,
DOI~10.1016/j.comgeo.2019.05.002

\bibitem[2021]{Messerschmidt21} The number of configurations of radii that
can occur in compact packings of the plane with discs of  sizes is finite.
arXiv:2110.15831v1 [math.MG] 29 Oct 2021.

\vskip2pt\hskip-1.5truecm{\sc{Moln\'ar, J.}}

\bibitem[1964]{Molnar64} Sui sistemi di punti con esigenza di spazio.
Atti Accad. Naz. Lincei Rend. Cl. Sci. Fis. Mat. Natur. (8) {\bf{36}} (1964) 336--339.
MR0170271 

\bibitem[1966a]{Molnar66a} Aggregati di cerchi di Minkowski.
Ann. Mat. Pura Appl. (4) {\bf{71}} (1966) 101--108. MR0205153, 
DOI~10.1007/BF02413737

\bibitem[1966b]{Molnar66b} Collocazioni di cerchi con esigenza di spazio.
Ann. Univ. Sci. Budapest. E\"otv\"os Sect. Math. {\bf{9}} (1966) 71--86. MR0205157 

\bibitem[1967b]{Molnar67b} On the $\lambda$-system of circles.
Acta Math. Acad. Sci. Hungar. {\bf{18}} (1967) 405--410. MR0218979, 
DOI~10.1007/BF02280300

\bibitem[1975]{Molnar75} On a generalisation of the Tammes problem.
Publ. Math. Debrecen {\bf{22}} (1975) no. 1-2, 109--114. MR0383259

\bibitem[1978]{Molnar78} Packing of congruent spheres in a strip.
Acta Math. Acad. Sci. Hungar. {\bf{31}} (1978) no. 1-2, 173--183. MR0487805,
DOI~10.1007/BF01896082

\vskip2pt\hskip-1.5truecm{\sc{Mughal, A.}} and {\sc{Weaire, D.}}

\bibitem[2014]{MughalWeaire} Theory of cylindrical dense packings of disks.
Phys. Rev. E {\bf{89}} no. 4 (2014) 042307.
10.1103/PhysRevE.89.042307

\vskip2pt\hskip-1.5truecm{\sc{Musin, O. R.}}

\bibitem[2019]{Musin19} Graphs and spherical two-distance sets.
European J. Combin. {\bf{80}} (2019) 311--325. MR3987372,
DOI~10.1016/j.ejc.2018.07.013

\vskip2pt\hskip-1.5truecm{\sc{Musin, O. R.}} and {\sc{Nikitenko, A. V.}}

\bibitem[2016]{MusinNikitenko} Optimal packings of congruent circles on a square flat torus.
Discrete Comput. Geom. {\bf{55}} (2016) no. 1, 1--20. MR3439258,
DOI~10.1007/s00454-015-9742-6

\vskip2pt\hskip-1.5truecm{\sc{Nasz\'odi, M.: Pach, J.}} and {\sc{Swanepoel, K.}}

\bibitem[2017]{NaszodiPachSwanepoel} Arrangements of homothets of convex bodies.
Mathematika {\bf{63}} (2017) no. 2, 696--710. MR3706603,
DOI~10.1112/S0025579317000122

\vskip2pt\hskip-1.5truecm{\sc{Nasz\'odi, M.}} and {\sc{Swanepoel, K.}}

\bibitem[2018]{NaszodiSwanepoel} Arrangements of homothets of convex bodies II.
Contrib. Discrete Math. {\bf{13}} (2018) no. 2, 116--123. MR3897229,
DOI~10.11575/cdm.v13i2.62732

\vskip2pt\hskip-1.5truecm{\sc{Natarajan B. K.}}

\bibitem[1988]{Natarajan} On planning assemblies.
Proceedings of the Fourth Annual ACM Symposium on Computational Geometry, 299--308, 1988.
DOI~10.1145/73393.73424

\vskip2pt\hskip-1.5truecm{\sc{Nielsen, F.}}

\bibitem[1965]{Nielsen} On the sum of the distances between $n$ points on a sphere. (Danish)
Nordisk Mat. Tidskr. {\bf{13}} (1965) no. 1-2, 45-50. Zbl 0132.17403

\vskip2pt\hskip-1.5truecm{\sc{Pach, J.}}

\bibitem[1977]{Pach77} On the permeability problem.
Studia Sci. Math. Hungar. {\bf{12}} (1977) no. 3-4, 419--424. MR0607097 

\vskip2pt\hskip-1.5truecm{\sc{Papadimitriou, Ch. H.}} and {\sc{Yannakakis, M.}}

\bibitem[1989]{PapadimitriouYannakakis89} Shortest paths without a map (extended abstract).
Automata, languages and programming (Stresa, 1989)  610--620,
Lecture Notes in Comput. Sci., 372, Springer, Berlin, 1989. MR1037079,
DOI~10.1007/BFb0035787

\bibitem[1991]{PapadimitriouYannakakis91} Shortest paths without a map.
16th International Colloquium on Automata, Languages, and Programming (Stresa, 1989).
Theoret. Comput. Sci. 84 (1991), no. 1, Algorithms
Automat. Complexity Games, 127--150. MR1122649,
DOI~10.1016/0304-3975(91)90263-2

\vskip2pt\hskip-1.5truecm{\sc{Peikert, R.}}

\bibitem[1994]{Peikert} Dichteste Packungen von gleichen Kreisen in einem Quad\-rat.
Elem. Math.  {\bf{49}}  (1994)  no. 1, 16--26. MR1261755

\vskip2pt\hskip-1.5truecm{\sc{Peikert, R.; W\"urtz, D.; Monagan, M.}} and de {\sc{Groot, C.}}

\bibitem[1992]{Peikert+} Packing circles in a square: a review and new results.
System modelling and optimization (Z\"urich, 1991)  45--54,
Lecture Notes in Control and Inform. Sci., 180, Springer, Berlin, 1992. MR1182322,
DOI~10.1007/BFb0113271

\vskip2pt\hskip-1.5truecm{\sc{Pillichshammer, F.}}

\bibitem[2000]{Pillichshammer} On the sum of squared distances in the Euclidean plane.
Arch. Math. (Basel) {\bf{74}} (2000) no. 6, 472--480. MR1753546,
DOI~10.1007/PL00000428

\vskip2pt\hskip-1.5truecm{\sc{Polyanskii, A.}}

\bibitem[2017]{Polyanskii17} Pairwise intersecting homothets of a convex body.
Discrete Math. {\bf{340}} (2017) no. 8, 1950--1956. MR3648222,
DOI~10.1016/j.disc.2017.04.002

\vskip2pt\hskip-1.5truecm{\sc{Przeworski, A.}}

\bibitem[2006]{Przeworski06b} Packing disks on a torus.
Discrete Comput. Geom. {\bf{35}} (2006) no. 1, 159--174. MR2183495,
DOI~10.1007/s00454-005-1198-7

\vskip2pt\hskip-1.5truecm{\sc{Rankin, R. A.}}

\bibitem[1955]{Rankin55} The closest packing of spherical caps in $n$ dimensions.
Proc. Glasgow Math. Assoc. {\bf{2}} (1955) 139--144. MR0074013, 
DOI~10.1017/S2040618500033219

\vskip2pt\hskip-1.5truecm{\sc{Rold\'an-Pensado, E.}}

\bibitem[2013]{Roldan-Pensado} Paths on the doubly covered region of a covering of the plane by discs.
Studia Sci. Math. Hungar. {\bf{50}} (2013) 465--469. MR3187828,
DOI~10.1556/SScMath.50.2013.4.1243

\vskip2pt\hskip-1.5truecm{\sc{Ry\v{s}kov, S. S.}} and {\sc{Baranovski\u{\i}, E. P.}}

\bibitem[1975]{RyskovBaranovskii75} Solution of the problem of the least dense
lattice covering of five-dimensional space by equal spheres. (Russian)
Dokl. Akad. Nauk SSSR {\bf{222}} (1975) no. 1, 39--42.
translation in: Soviet Math. Dokl. {\bf{16}} (1975) 586--590. MR0427238

\bibitem[1976]{RyskovBaranovskii76} $C$-types of $n$-dimensional lattices and 5-dimensional
primitive parallelohedra (with application to the theory of coverings).
Cover to cover translation of Trudy Mat. Inst. Steklov 137 (1976).
Translated by R. M. Erdahl. Proc. Steklov Inst. Math. 1978, no. 4, 140 pp. MR0535314

\vskip2pt\hskip-1.5truecm{\sc{Schaer, J.}}

\bibitem[1966a]{Schaer66a} On the densest packing of spheres in a cube.
Canad. Math. Bull. {\bf{9}} (1966) 265--270. MR0200797, 
DOI~10.4153/CMB-1966-033-0

\bibitem[1966b]{Schaer66b} The densest packing of five spheres in a cube.
Canad. Math. Bull. {\bf{9}} (1966) 271--274. MR0200798, 
DOI~10.4153/CMB-1966-034-8

\bibitem[1966c]{Schaer66c} The densest packing of six spheres in a cube.
Canad. Math. Bull. {\bf{9}} (1966) 275--280. MR0202055, 
DOI~10.4153/CMB-1966-035-5

\bibitem[1994]{Schaer94} The densest packing of ten congruent spheres in a cube.
Intuitive geometry (Szeged, 1991) 403--424, Colloq. Math. Soc. J\'anos Bolyai, 63,
North-Holland, Amsterdam, 1994. MR1383635 

\vskip2pt\hskip-1.5truecm{\sc{Sch\"urmann, A.}}

\bibitem[2002a]{Schurmann02a} On extremal finite packings.
Discrete Comput. Geom. {\bf{28}} (2002) no. 3, 389--403. MR1923959,
DOI~10.1007/s00454-002-0747-6

\vskip2pt\hskip-1.5truecm{\sc{Sch\"urmann, A.}} and {\sc{Vallentin, F.}}

\bibitem[2006]{SchurmannVallentin} Computational approaches to lattice packing and covering problems.
Discrete Comput. Geom. {\bf{35}} (2006) no. 1, 73--116. MR2183491,
DOI~10.1007/s00454-005-1202-2

\vskip2pt\hskip-1.5truecm{\sc{Shephard, G. C.}}

\bibitem[1970]{Shephard70} On a problem of Fejes T\'oth.
Studia Sci. Math. Hungar. {\bf{5}} (1970) 471--473. MR0285995 

\vskip2pt\hskip-1.5truecm{\sc{Snoeyink, J.}} and {\sc{Stolfi, J.}}

\bibitem[1994]{SnoeyinkStolfi} Objects that cannot be taken apart with two hands.
ACM Symposium on Computational Geometry (San Diego, CA, 1993).
Discrete Comput. Geom. {\bf{12}} (1994) no. 3, 367--384. MR1298917, 
DOI~10.1007/BF02574386

\vskip2pt\hskip-1.5truecm{\sc{Sperling, G.}}

\bibitem[1960]{Sperling} L\"osung einer elementargeometrischen Frage von Fejes T\'oth.
Arch. Math. (Basel) {\bf{11}} (1960) 69--71. MR0112077, 
DOI~10.1007/BF01236910

\vskip2pt\hskip-1.5truecm{\sc{Stolarsky, K. B.}}

\bibitem[1972]{Stolarsky72} Sums of distances between points on a sphere.
Proc. Amer. Math. Soc. {\bf{35}} (1972) 547--549. MR0303418,
DOI~10.1090/S0002-9939-1972-0303418-3

\bibitem[1973]{Stolarsky73} Sums of distances between points on a sphere. II.
Proc. Amer. Math. Soc. {\bf{41}} (1973) 575--582. MR0333995,
DOI~10.1090/S0002-9939-1973-0333995-9

\vskip2pt\hskip-1.5truecm{\sc{Szab\'o, P. G.}}; {\sc{Mark\'ot, M. Cs.}}; {\sc{Csendes, T.}}; {\sc{Specht, E.}}; {\sc{Casado, L. G.}} and
{\sc{Garc\'{\i}a, I.}}

\bibitem[2007]{Szabo+} New approaches to circle packing in a square. With program codes. With 1 CD-ROM (Linux).
Springer Optimization and Its Applications, 6. Springer, New York, 2007. xiv+238 pp. MR2296052 

\vskip2pt\hskip-1.5truecm{\sc{Szab\'o, L.}} and {\sc{Ujv\'ary-Menyh\'art, Z.}}

\bibitem[2002]{SzaboUjvary-Menyhart} Clouds of planar convex bodies.
Aequationes Math. {\bf{63}} (2002) no. 3, 292--302. MR1904721,
DOI~10.1007/s00010-002-8025-5

\vskip2pt\hskip-1.5truecm{\sc{Talata, I.}}

\bibitem[2000a]{Talata00a} On translational clouds for a convex body.
Geom. Dedicata  {\bf{80}}  (2000)  no. 1-3, 319--329. MR1762518, 
DOI~10.1023/A:1005279901749

\vskip2pt\hskip-1.5truecm{\sc{Torquato, S.}} and {\sc{Jiao, Y.}}

\bibitem[2009a]{TorquatoJiao09a} Dense packings of the Platonic and Archimedean solids.
Nature {\bf{460}} (2009) 876--879.
DOI~10.1038/nature08239

\bibitem[2009b]{TorquatoJiao09b} Dense packings of polyhedra: Platonic and Archimedean solids.
Phys. Rev. E (3) {\bf{80}} (2009) no. 4, 041104, 21 pp. MR2607418,
DOI~10.1103/PhysRevE.80.041104

\bibitem[2009c]{TorquatoJiao09c} Analytical constructions of a family
of dense tetrahedron packings and the role of symmetry.
arXiv:0912.4210v1 [cond-mat.stat-mech] 21 Dec 2009

\vskip2pt\hskip-1.5truecm{\sc{V\'as\'arhelyi, \'E.}}

\bibitem[2003]{Vasarhelyi03} Experimentieren um einen Satz zu finden --- vollst\"anding separierbare Mosaike auf der Kugel
und ihre Anwendungen, Teaching Math. and Comput. Sci. {\bf{1}} (2003) no. 2, 297--319.

\vskip2pt\hskip-1.5truecm{\sc{Vermes, I.}}

\bibitem[1996]{Vermes96} Total-separable Kreissysteme und
Mosaike in der hyperbolischen Ebene.
Studia Sci. Math. Hungar. {\bf{32}} (1996) no. 3-4, 377--382. MR1432179

\vskip2pt\hskip-1.5truecm{\sc{Witsenhausen, H. S.}}

\bibitem[1974]{Witsenhausen} On the maximum of the sum of squared distances under a diameter constraint.
Amer. Math. Monthly {\bf{81}} (1974) 1100--1101. MR0370367,
DOI~10.1080/00029890.1974.11993743

\vskip2pt\hskip-1.5truecm{\sc{Xue, F.}} and {\sc{Zong, C.}}

\bibitem[2018]{XueZong} On lattice coverings by simplices.
Adv. Geom. {\bf{18}} (2018) no. 2, 181--186. MR3785419,
DOI~10.1515/advgeom-2017-0049

\vskip2pt\hskip-1.5truecm{\sc{Ziegler, G. M.}}

\bibitem[2010]{Ziegler10} Three Mathematics Competitions.
In An Invitation to Mathematics From Competitions to Research 195--205.
Dierk Schleicher and Malte Lackmann Editors,
Springer Heidelberg Dordrecht London New York, 2010. Zbl 1323.00022,
DOI~10.1007/978-3-642-19533-4\_13

\vskip2pt\hskip-1.5truecm{\sc{Zong, C.}}

\bibitem[1997b]{Zong97b} A problem of blocking light rays.
Geom. Dedicata {\bf{67}} (1997) no. 2, 117--128. MR1474113,
DOI~10.1023/A:1004974630549

\bibitem[1999a]{Zong99a} A note on Hornich's problem.
Arch. Math. (Basel) {\bf{72}} (1999) no. 2, 127--131. MR1667053, 
DOI~10.1007/s000130050313

\bibitem[2003]{Zong03} Simultaneous packing and covering in three-dimensional Euclidean space.
J. London Math. Soc. (2)  {\bf{67}}  (2003)  no. 1, 29--40. MR1942409, 
DOI~10.1112/S0024610702003873

\bibitem[2008]{Zong08} The simultaneous packing and covering constants in the plane.
Adv. Math. {\bf{218}} (2008) no. 3, 653--672. MR2414317,
DOI~10.1016/j.aim.2008.01.007

\bibitem[2019]{Zong19} A computer approach to determine the densest translative tetrahedron packings.
Experimental Mathematics (2019)
DOI~10.1080/10586458.2019.1582378


\end{thebibliography}
\end{document}
